\documentclass{sbc}
\smartqed  % flush right qed marks, e.g. at end of proof
\usepackage{graphicx}
\usepackage{mathptmx}      % use Times fonts if available on your TeX system
\usepackage{latexsym}
\usepackage{amsmath} 
\usepackage{amssymb,psfrag}  
\newcommand{\vp}{\tau}
\newcommand{\emvp}{e^{-\tau}}%{\hbox{$e^{-\tau\hskip-4pt\iota}$}}
\newcommand{\havp}{\hat\vp}
\newcommand{\haro}{\hat\rho}
\newcommand{\navp}{\nabla\hskip-1pt\vp}
\newcommand{\df}{d\hskip-.8pt\fx}
\newcommand{\fc}{H}
\newcommand{\bbI}{\mathbf{I}}
\newcommand{\bbR}{\mathrm{I\!R}}
\newcommand{\bbRP}{\bbR\mathrm{P}}
\newcommand{\bbC}{{\mathchoice {\setbox0=\hbox{$\displaystyle\mathrm{C}$}
\hbox{\hbox to0pt{\kern0.4\wd0\vrule height0.9\ht0\hss}\box0}} 
{\setbox0=\hbox{$\textstyle\mathrm{C}$}\hbox{\hbox 
to0pt{\kern0.4\wd0\vrule height0.9\ht0\hss}\box0}} 
{\setbox0=\hbox{$\scriptstyle\mathrm{C}$}\hbox{\hbox 
to0pt{\kern0.4\wd0\vrule height0.9\ht0\hss}\box0}} 
{\setbox0=\hbox{$\scriptscriptstyle\mathrm{C}$}\hbox{\hbox 
to0pt{\kern0.4\wd0\vrule height0.9\ht0\hss}\box0}}}} 
\newcommand{\bbCP}{\bbC\mathrm{P}}

\newcommand{\dimc}{\dim_{\hskip.4pt\bbC\hskip-1.2pt}^{\phantom i}}
\newcommand{\Lie}{\pounds}
\newcommand{\bz}{b\hs}
\newcommand{\lj}{\langle}
\newcommand{\rg}{\rangle}
\newcommand{\lr}{\lj\hskip1.6pt,\rg}
\newcommand{\cj}{c}
\newcommand{\fx}{f}
\newcommand{\hs}{\hskip.7pt}
\newcommand{\hh}{\hskip.4pt}
\newcommand{\hn}{\hskip-.4pt}
\newcommand{\nh}{\hskip-.7pt}
\newcommand{\nnh}{\hskip-1.5pt}
\newcommand{\tm}{{T\hskip-.3ptM}}

\newcommand{\yj}{c}
\newcommand{\kfo}{\omega}
\newcommand{\hatg}{{\hat{g\hskip2pt}\hskip-1.3pt}}
\newcommand{\hato}{{\hat{\kfo}}}
\newcommand{\hasc}{{\hat{\se}}}
\newcommand{\hanb}{\hat\nabla}

\newcommand{\sa}{\sigma}

\newcommand{\fy}{{\theta}}

\newcommand{\omh}{\omega^{(h)}}
\newcommand{\ve}{\varepsilon}

\newcommand{\se}{\mathrm{s}}
\newcommand{\proj}{\pi}

\begin{document}

\renewcommand{\theequation}{\thesection.\arabic{equation}}

\voffset=32pt\hoffset=53.4pt  %for arXiv & homepage

\title{Special bi\-con\-for\-mal changes of K\"ah\-ler surface metrics}

\titlerunning{Bi\-con\-for\-mal changes of K\"ah\-ler metrics}

\author{Andrzej Derdzinski}

\institute{A. Derdzinski \at
Department of Mathematics, The Ohio State University, Columbus, OH 43210, USA\\
Tel.: +1-614-2924012\\
Fax: +1-614-2921479\\
\email{andrzej@math.ohio-state.edu}}

\maketitle

\begin{abstract}
The term ``special bi\-con\-for\-mal change'' refers, 
basically, to the situation where a given nontrivial real-hol\-o\-mor\-phic 
vector field on a complex manifold is a gradient relative to two K\"ah\-ler 
metrics, and, simultaneously, an eigen\-vec\-tor of one of the metrics 
treated, with the aid of the other, as an en\-do\-mor\-phism of the tangent 
bundle. A special bi\-con\-for\-mal change is called nontrivial if the two 
metrics are not each other's constant multiples. For instance, according to 
a 
1995 result of LeBrun, a nontrivial special bi\-con\-for\-mal change exists for the 
con\-for\-mal\-ly-Ein\-stein K\"ah\-ler metric on the two-point 
blow-up of the complex projective plane, recently discovered by Chen, LeBrun 
and Weber; the real-hol\-o\-mor\-phic vector field involved is the gradient of 
its scalar curvature. The present paper establishes the existence of 
nontrivial special bi\-con\-for\-mal changes for some canonical metrics on Del 
Pezzo surfaces, viz.\ K\"ahler-Einstein metrics (when a nontrivial 
hol\-o\-mor\-phic vector field exists), non-Ein\-stein K\"ah\-ler-Ric\-ci 
sol\-i\-tons, and K\"ah\-ler metrics admitting %with 
nonconstant Kil\-ling potentials %that have 
with geodesic gradients.% An appendix provides a classification of compact 
%K\"ah\-ler surfaces admitting Kil\-ling potentials with this last property.
\end{abstract}

\setcounter{section}{0}
\section{Introduction}%\label{in}
\setcounter{equation}{0}
By a {\it met\-ric-\hh po\-ten\-tial pair\/} on a complex manifold $\,M\,$ 
with $\,\dimc M\ge2\,$ we mean any pair $\,(g,\nnh\vp)\,$ formed by a 
K\"ah\-ler metric $\,g\,$ on $\,M\,$ and a nonconstant Kil\-ling potential 
$\,\vp\,$ for $\,g$, that is, a function $\,\vp:M\to\bbR\,$ such that 
$\,J(\navp)\,$ is a nontrivial Kil\-ling field on the K\"ah\-ler manifold 
$\,(M,g)$. Another met\-\hbox{ric-}\hskip.7ptpo\-ten\-tial pair $\,(\hatg,\havp)\,$ on the same 
complex manifold $\,M\,$ is said to arise from $\,(g,\nnh\vp)\,$ by a {\it 
special bi\-con\-for\-mal change\/} if
\begin{equation}\label{bcm}
\mathrm{i)}\hskip12pt
\hatg\,\,=\,\,\fx g\,\,-\,\,\fy\hs(d\vp\otimes d\vp\hs+\hs\xi\otimes\xi)\hs,
\hskip16pt\mathrm{ii)}\hskip12pt\hanb\havp\,=\,\nabla\vp
\end{equation}
for $\,\xi=g(J(\navp),\,\cdot\,)\,$ and some $\,C^\infty$ functions 
$\,\fx,\fy:M\to\bbR$. The equality in (\ref{bcm}.ii) states that the 
$\,\hatg$-gra\-di\-ent of $\,\havp\,$ coincides with the $\,g$-gra\-di\-ent of 
$\,\vp$.

A special bi\-con\-for\-mal change as above will be called {\it trivial\/} if 
$\,\fx\,$ is a positive constant, $\,\fy=\hs0$, and $\,\havp\,$ equals 
$\,\fx\vp\,$ plus a constant.

Gan\-chev and Mi\-ho\-va \cite[Section~4]{ganchev-mihova} studied 
bi\-con\-for\-mal changes of a more general type. In their approach, 
$\,\vp:M\to\bbR\,$ is not required to be a Kil\-ling potential.

The existence of nontrivial special bi\-con\-for\-mal changes has already been 
established for some met\-\hbox{ric-}\hskip.7ptpo\-ten\-tial pairs $\,(g,\nnh\vp)$. LeBrun 
\cite{lebrun} proved it when $\,g\,$ is a K\"ah\-ler metric on a compact 
complex surface, con\-for\-mal to a non-K\"ah\-ler Ein\-stein metric, and 
$\,\vp\,$ is the scalar curvature of $\,g$. Both the one-point and two-point 
blow-ups of $\,\bbCP^2$ are known to admit metrics with the properties just 
listed (the latter, due to a recent result of Chen, LeBrun and Weber 
\cite{chen-lebrun-weber}; see also Section~\ref{ce}). On the other hand, 
Gan\-chev and Mi\-ho\-va \cite{ganchev-mihova} exhibited a nontrivial special 
bi\-con\-for\-mal change leading from $\,(g,\nnh\vp)$, for any non\-flat 
K\"ah\-ler metric $\,g\,$ of qua\-si-con\-stant hol\-o\-mor\-phic sectional 
curvature, and suitable $\,\vp$, to a met\-\hbox{ric-}\hskip.7ptpo\-ten\-tial pair 
$\,(\hatg,\havp)\,$ in which the K\"ah\-ler metric $\,\hatg\,$ is flat.

This paper addresses the existence question for nontrivial special 
bi\-con\-for\-mal changes of met\-\hbox{ric-}\hskip.7ptpo\-ten\-tial pairs in complex 
dimension $\,2$. It is not known whether all met\-\hbox{ric-}\hskip.7ptpo\-ten\-tial pairs 
$\,(g,\nnh\vp)\,$ on compact complex surfaces admit such changes. However, 
nontrivial special bi\-con\-for\-mal changes of $\,(g,\nnh\vp)\,$ always exist 
locally, at points where $\,d\vp\ne0\,$ (Remark~\ref{exloc} at the end of 
Section~\ref{sb}).

Bi\-con\-for\-mal changes of a more general kind than those defined above  
are introduced in Section~\ref{os}, where it is also shown that such a 
generalized bi\-con\-for\-mal change exists between any two 
$\,\mathrm{U}(2)$-in\-var\-i\-ant K\"ah\-ler metrics on $\,\bbCP^2$ or on 
the one-point blow-up of $\,\bbCP^2\nnh$.

Theorems~\ref{dftef} and~\ref{opidd}, stated and proved in Sections~\ref{og} 
and~\ref{ac}, provide two general mechanisms allowing one to construct 
examples of nontrivial special bi\-con\-for\-mal changes. They are based on 
criteria for the existence of such changes that are, in addition, required to 
satisfy a certain functional dependence relation, or to yield a metric in the 
same K\"ah\-ler class; in the former case the criterion amounts to a Laplacian 
condition. 

The first main result of the paper, derived from Theorem~\ref{dftef}, is the 
existence of nontrivial special bi\-con\-for\-mal changes of various 
canonical metrics on Del Pezzo surfaces. Specifically, they are shown to exist 
for all met\-\hbox{ric-}\hskip.7ptpo\-ten\-tial pairs $\,(g,\nnh\vp)\,$ with suitably chosen 
$\,\vp$, on compact complex surfaces $\,M\nh$, such that $\,g\,$ is
\begin{enumerate}
  \def\theenumi{{\rm\roman{enumi}}}
\item[(i)] any K\"ah\-ler-Ein\-stein metric with positive scalar curvature 
(and $\,M\,$ admits a nontrivial hol\-o\-mor\-phic vector field), or
\item[(ii)] any non-Ein\-stein K\"ah\-ler-Ric\-ci soliton, or
%\item[$\bullet$] K\"ah\-ler metric nontrivially con\-for\-mal to an 
%Ein\-stein metric, or
\item[(iii)] any K\"ah\-ler metric admitting a special K\"ah\-ler-Ric\-ci 
potential $\,\vp$.%, as well as for
\end{enumerate}
The second main result is Theorem~\ref{geodg}, establishing the existence of 
nontrivial special bi\-con\-for\-mal changes of $\,(g,\nnh\vp)\,$ whenever 
$\,(M,g)\,$ is a compact K\"ah\-ler surface and the integral curves of 
$\,\navp\,$ are re\-pa\-ram\-e\-trized geodesics. Being a special 
K\"ah\-ler-Ric\-ci potential is sufficient for $\,\vp\,$ to have this last 
property, but it is not necessary; more general examples are described in the 
Appendix.

Two K\"ah\-ler metrics on a given complex surface cannot be nontrivially 
con\-for\-mal. The relation of ``general bi\-con\-for\-mal equivalence'' is 
not of much interest here either, since it holds locally, almost everywhere, 
for any two K\"ah\-ler surface metrics (Section~\ref{sb}). On the other hand, 
on compact complex surfaces, a special bi\-con\-for\-mal change between two 
given met\-\hbox{ric-}\hskip.7ptpo\-ten\-tial pairs exists sometimes, though not very often, 
and if it does exist, it amounts to an explicit description of one K\"ah\-ler 
metric in terms of the other. For instance, as shown at the end of 
Section~\ref{ut}, the one-point blow-up of $\,\bbCP^2$ admits a 
bi\-con\-for\-mal change of a more general type, introduced in 
Section~\ref{os}, leading from the K\"ah\-ler-Ric\-ci sol\-i\-ton constructed 
by Koiso \cite{koiso} and, independently, Cao \cite{cao}, to one of 
Ca\-la\-bi's extremal K\"ah\-ler metrics \cite{calabi}, con\-for\-mal to the 
non-K\"ah\-ler, Ein\-stein metric found by Page \cite{page}.

\section{Preliminaries}\label{pr}
\setcounter{equation}{0}
All manifolds and Riemannian metrics are assumed to be of class 
$\,C^\infty\nnh$. A manifold is by definition connected.

Given a Riemannian manifold $\,(M,g)$, the {\it divergence\/} of a vector 
field $\,w\,$ or a bundle morphism $\,A:\tm\to\tm\,$ is defined  as usual, by 
$\,\mathrm{div}\hskip2ptw=\hs\mathrm{tr}\,\nabla\nh w\,$ and 
$\,\mathrm{div}\hskip2ptA=\xi$, for the $\,1$-form $\,\xi\,$ sending any 
vector field $\,w\,$ to the function 
$\,\xi(w)=\mathrm{div}\,(\hn Aw)-\hs\mathrm{tr}\,(\hn A\hs\nabla\nh w)$. The 
inner product $\,\lr\,$ of $\,2$-forms is characterized by 
$\,2\lj\sa,\sa\rg=-\hs\mathrm{tr}\,A^2\nnh$, where $\,A:\tm\to\tm\,$ is the 
bundle morphism with $\,g(\hn Aw,\,\cdot\,)=\sa(w,\,\cdot\,)\,$ for all vector 
fields $\,w$. In coordinates, $\,\mathrm{div}\hskip2ptw=w^{\hs j}_{,\hs j}$, 
$\,(\mathrm{div}\hskip2ptA)_j\nh=A^{\nh k}_{j,\hs k}$ and 
$\,2\lj\sa,\sa\rg=\hs\sa_{jk}\hs\sa\hs^{jk}\nnh$. Also, for any 
$\,2$-form $\,\sa\,$ and vector fields $w,w\hh'\nnh$,
\begin{equation}\label{oab}
\lj\sa,\alpha\wedge\alpha\hh'\hh\rg\,=\,\sa(w,w\hh'\hh)\hs,\hskip13pt
\mathrm{where}\hskip6pt\alpha=g(w,\,\cdot\,),\hskip4pt\alpha\hh'\nh
=g(w\hh'\nnh,\,\cdot\,)\hh.
\end{equation}
\begin{lemma}\label{cinff}Suppose that\/ $\,\delta,\ve\in(0,\infty)\,$ and\/ 
$\,\vp,\psi:(-\hh\delta,\ve)\to\bbR\,$ are\/ $\,C^\infty$ functions such that, 
if the dot stands for the derivative with respect to the variable\/ 
$\,t\in(-\hh\ve,\ve)$,
\begin{enumerate}
  \def\theenumi{{\rm\alph{enumi}}}
\item $\dot\vp(0)=0\ne\ddot\vp(0)\,$ and\/ $\,\dot\vp\hs\ddot\vp\ne0\,$ 
everywhere in\/ $\,(-\hh\delta,0)\cup(0,\ve)$,
\item $\,\vp:(-\hh\delta,0)\to\bbR\,$ and\/ 
$\,\vp:(0,\ve)\to\bbR\,$ both have the same range\/ $\,\bbI\subset\bbR$,
\item $\psi(t)=G(\vp(t))\,$ for some function\/ $\,G:\bbI\to\bbR\,$ and all\/ 
$\,t\in(-\hh\delta,0)\cup(0,\ve)$.
\end{enumerate}
Then\/ $\,G\,$ has a\/ $\,C^\infty$ extension to the half-open interval\/ 
$\,\bbI\hs\cup\vp(0)$.
\end{lemma}
\begin{proof}One can view $\,\vp\,$ as a new $\,C^\infty$ coordinate on both 
$\,(-\hh\ve,0)\,$ and $\,(0,\ve)$. Thus, $\,G:\bbI\to\bbR\,$ is of class 
$\,C^\infty$, and so are all the derivatives $\,d\hh^k\nh G\hs/d\vp^k$ 
treated as functions on $\,\bbI$. Let us prove by induction on $\,k\ge0\,$ that 
$\,d\hh^k\nh G\hs/d\vp^k$ is a $\,C^\infty$ function of the variable 
$\,t\in(-\hh\ve,\ve)\,$ (and, in particular, has a limit at the endpoint 
$\,\vp(0)\,$ of $\,\bbI$). The induction step: by (b) -- (c), 
$\,\chi=d\hh^k\nh G\hs/d\vp^k\,$ treated as a $\,C^\infty$ function on 
$\,(-\hh\delta,\ve)\,$ has the same range on $\,(-\hh\delta,0)\,$ as on 
$\,(0,\ve)$, which also remains true when $\,\delta,\ve\,$ are replaced with 
suitably related smaller positive numbers $\,\delta\hh'\nnh,\ve\hh'\nnh$, and 
such $\,\delta\hh'\nnh,\ve\hh'$ may be chosen arbitrarily close to $\,0$. 
Hence $\,\dot\chi(0)=0$. As $\,\dot\vp\,$ is a new $\,C^\infty$ coordinate on 
$\,(-\hh\ve,\ve)$, vanishing at $\,0$, smooth functions on $\,(-\hh\ve,\ve)\,$ 
that vanish at $\,0\,$ are smoothly divisible by $\,\dot\vp$. Consequently, 
$\,d\hh^{k+1}\nh G\hs/d\vp^{\hs k+1}=d\chi\hs/d\vp=\dot\chi\hs/\dot\vp\,$ is 
a smooth function of $\,t\in(-\hh\ve,\ve)$.
\qed
\end{proof}
\begin{remark}\label{divis}Let $\,F\,$ be a $\,C^\infty$ function 
$\,\,U\nh\times D\to\bbR$, where $\,\,U\nh\subset\bbR^k$ is an open set and 
$\,D\subset\bbC\,$ is a disk centered at $\,0$. If $\,F(y,0)=0\,$ and 
$\,F(y,zq)=F(y,z)\,$ for all $\,(y,z,q)\in\bbR^k\nnh\times\bbC^2$ with 
$\,|q|=1$, then $\,F(y,z)=|z|^2h(y,z)\,$ for some $\,C^\infty$ function 
$\,h:U\nh\times D\to\bbR\hh$. If, in addition, the Hess\-i\-an of $\,F\,$ is 
nonzero everywhere in $\,\,U\nh\times\{0\}$, then so is $\,h$.

In fact, for $\,r\in\bbR\,$ close to $\,0$, the function 
$\,(y,r)\mapsto F(y,r)\,$ is smooth and vanishes when $\,r=0$, so that it 
is smoothly divisible by $\,r\,$ (due to the first-or\-der Taylor formula). 
The same applies to $\,(y,r)\mapsto F(y,r)/r$. The last claim holds since, on 
$\,\,U\nh\times\{0\}$, the Hess\-i\-an of $\,F\,$ equals $\,2\hh h\,$ times 
the Euclidean metric of $\,\bbC$.
\end{remark}
We will use the {\it connectivity lemma\/} for \hbox{Morse\hs-}\hskip0ptBott 
functions $\,\vp\,$ on compact manifolds $\,M$, stating that, if the positive 
and negative indices of the Hess\-i\-an of $\,\vp\,$ at every critical point 
are both different from $\,1$, then the $\,\vp$-pre\-im\-age of every real 
number is connected. See \cite[Lemma 3.46 on p.\ 124]{nicolaescu}.

\section{K\"ah\-ler manifolds}\label{km}
\setcounter{equation}{0}
Let $\,M\,$ be a complex manifold. Its com\-plex-struc\-ture tensor is always 
denoted by $\,J$. Given a real $\,1$-form $\,\mu\,$ on $\,M\nh$, the symbol 
$\,\mu J\,$ stands for the $\,1$-form $\,J^*\nnh\mu$, so that
\begin{equation}\label{mux}
(\mu J)_x\hs=\,\,\mu_x\circ J_x:T_xM\to\bbR\
\end{equation}
at any point $\,x\in M\nh$. If $\,g\,$ is a K\"ah\-ler metric on $\,M\nh$, the 
K\"ah\-ler form of $\,g\,$ 
is $\,\kfo=g(J\,\cdot\,,\,\cdot\,)$, while $\,\nabla\,$ denotes both the 
Le\-vi-Ci\-vi\-ta connection of $\,g\,$ and the $\,g$-gra\-di\-ent. 
Real-hol\-o\-mor\-phic vector fields on $\,M\,$ are the sections $\,v\,$ of 
$\,\tm\,$ such that $\,\Lie_v\hs J=0$, which, for any fixed K\"ah\-ler 
metric $\,g\,$ on $\,M\nh$, is equivalent to $\,[\hh J,\nabla v\hh]=0$. The 
commutator $\,[\hskip2.2pt,\hskip1pt]\,$ is applied here to $\,J\,$ and 
$\,\nabla v\,$ treated as vec\-tor-bun\-dle morphisms $\,\tm\to\tm$, the 
latter acting on vector fields $\,w\,$ by 
$\,(\nabla v)w=\nabla_{\hskip-1.2ptw}v$. The fact that
\begin{equation}\label{idd}
2\hs i\hskip1pt\partial\overline{\partial}\hskip0pt\psi\,
=\,(\nabla\nh d\psi)(J\,\cdot\,,\,\cdot\,)\,
-\,(\nabla\nh d\psi)(\,\cdot\,,J\,\cdot\,)
\end{equation}
for any K\"ah\-ler metric $\,g\,$ on $\,M\nh$, and any $\,C^2$ function 
$\,\psi:M\to\bbR$, will be used below and in Section~\ref{ac}. In the 
following (well-known) lemmas, 
$\,\imath_v\hh\bz=\bz(v,\,\cdot\,,\dots,\,\cdot\,)\,$ for vector fields 
$\,v\,$ and co\-var\-i\-ant tensor fields $\,\bz$, while 
$\,[\hh d(d_v\psi)]J\,$ is defined by (\ref{mux}) with $\,\mu=d(d_v\psi)$.
\begin{lemma}\label{ivddf}In a K\"ah\-ler manifold\/ $\,(M,g)\,$ one has\/ 
$\,2\hh\imath_v(i\hskip1pt\partial\overline{\partial}\hs\psi)\hs
=\hs d(d_{Jv}\psi)-\hs[\hh d(d_v\psi)]J\,$ for any real-hol\-o\-mor\-phic 
vector field\/ $\,v\,$ on\/ $\,M\,$ and any\/ $\,C^2$ function\/ 
$\,\psi:M\to\bbR\hh$. 
\end{lemma}
\begin{proof}Let us set $\,u=Jv\,$ and $\,\mu=
2\hh\imath_v(i\hskip1pt\partial\overline{\partial}\hs\psi)$. 
By (\ref{idd}), 
$\,\mu J=[\imath_u(\nabla d\psi)]J+\imath_v(\nabla d\psi)$. In coordinates, 
this reads 
$\,(\mu J)_k^{\phantom i}=\psi_{,pq}v^{\hh s}J_s^pJ_k^q
+\psi_{,sk}^{\phantom i}v^{\hh s}\nnh$. However,
$\,\psi_{,pq}^{\phantom i}v^{\hh s}J_s^p
=(\psi_{,p}^{\phantom i}v^{\hh s})_{,q}^{\phantom i}J_s^p
-\psi_{,p}^{\phantom i}v^{\hh s}_{,q}J_s^p$, while, as 
$\,v\,$ is real-hol\-o\-mor\-phic, $\,v^{\hh s}_{,q}J_s^p=J_q^sv^p_{,s}$. 
Hence $\,(\mu J)_k^{\phantom i}
=(\psi_{,p}^{\phantom i}v^{\hh s}J_s^p)_{,q}^{\phantom i}J_k^q
+\psi_{,s}^{\phantom i}v^{\hh s}_{,k}+\psi_{,sk}^{\phantom i}v^{\hh s}\nnh
=\{[d(d_u\psi)]J+d(d_v\psi)\}_k^{\phantom i}$.
\qed
\end{proof}
\begin{lemma}\label{hnteq}Let there be given\/ $\,C^\infty$ functions\/ 
$\,\vp,\psi:M\to\bbR\,$ on a complex manifold\/ $\,(M,g)\,$ and two 
K\"ah\-ler metrics $\,g,\hatg\,$ on\/ $\,M\,$ such that the K\"ah\-ler forms\/ 
$\,\kfo\,$ of\/ $\,g\,$ and\/ $\,\hato\,$ of\/ $\,\hatg\,$ are related 
by\/ $\,\hato=\kfo\hs+2i\hskip1pt\partial\overline{\partial}\hs\psi$. If the 
the\/ $\,g$-gra\-di\-ent\/ $\,v\,$ of\/ $\,\vp\,$ is real-hol\-o\-mor\-phic 
and\/ $\,d_{Jv}\psi=0$, then\/ $\,v\,$ is also the\/ $\,\hatg$-gra\-di\-ent 
of\/ $\,\havp=\vp+d_v\psi$.
\end{lemma}
\begin{proof}By (\ref{mux}), 
$\,(\imath_v\kfo)J=\kfo(v,J\,\cdot\,)=g(Jv,J\,\cdot\,)=g(v,\,\cdot\,)=d\vp$. 
Consequently, $\,(\imath_v\hato)J=\hatg(v,\,\cdot\,)$. Lemma~\ref{ivddf} now 
yields $\,\hatg(v,\,\cdot\,)=(\imath_v\hato)J=d\hs(\vp+d_v\psi)$.
\qed
\end{proof}
Here is another well-known lemma.
\begin{lemma}\label{cltwf}A differentiable $\,2$-form\/ $\,\eta\,$ on a 
K\"ah\-ler surface\/ $\,(M,g)\,$ is closed if and only if\/ 
$\,d\hh\lj\hs\omega,\eta\rg=-\hs\mathrm{div}\hs JA$, where\/ $\,J\,$ is the 
complex structure, $\,\omega\,$ denotes the K\"ah\-ler form of\/ $\,g$, and\/ 
$\,A:\tm\to\tm\,$ is the bundle morphism with 
$\,g(\hn Av,\,\cdot\,)=\eta(v,\,\cdot\,)\,$ for all vector fields\/ $\,v$.
\end{lemma}
\begin{proof}The operator $\,\imath_\omega$ sending every differential 
$\,3$-form $\,\zeta\,$ on $\,M\,$ to the $\,1$-form 
$\,\imath_\omega\hs\zeta\,$ such that $\,(\imath_\omega\hs\zeta)(v)
=\lj\hs\omega,\hs\zeta(v,\,\cdot\,,\,\cdot\,)\rg\,$ for all vector fields 
$\,v\,$ is, by dimensional reasons, an isomorphism, since 
$\,\imath_\omega(\xi\wedge\omega)=\xi\,$ for any $\,1$-form $\,\xi$. The 
assertion now follows from the lo\-cal-co\-or\-di\-nate formula 
$\,2\hh(\imath_\omega\hs d\hh\eta)_j\nh
=\omega\hh^{kl}(\eta_{kl,\hs j}\nh+\eta_{lj,\hs k}\nh+\eta_{jk,\hs l})$.
\qed
\end{proof}

\section{Kil\-ling potentials}\label{kp}
\setcounter{equation}{0}
Let $\,\vp\,$ be a Kil\-ling potential on a compact K\"ah\-ler manifold 
$\,(M,g)$. As usual, this means that $\,\vp\,$ is a $\,C^\infty$ function 
$\,M\to\bbR\,$ and $\,J(\navp)\,$ is a Kil\-ling field on $\,(M,g)$. In 
other words, $\,\navp\,$ is a real-hol\-o\-mor\-phic vector field, or, 
equivalently, the $\,2$-ten\-sor field $\,\nabla d\vp\,$ is 
Her\-mit\-i\-an. Using the notation
\begin{equation}\label{not}
v=\navp\hs,\hskip17ptu=Jv\hs,\hskip17pt\xi=g(u,\,\cdot\,)\hs,\hskip17ptQ
=g(v,v)\hs,\hskip17ptY\nnh=\Delta\vp\hs,
\end{equation}
here and and throughout the paper, one then has
\begin{equation}\label{rcv}
\mathrm{a)}\hskip6pt2\hskip1.5pt\mathrm{Ric}\hh(v,\,\cdot\,)=-\hs d\hskip.2ptY,\hskip18pt
\mathrm{b)}\hskip6pt2\hh\nabla d\vp(v,\,\cdot\,)=dQ\hh.
\end{equation}
In fact, (\ref{rcv}.b) holds for any $\,C^\infty$ function $\,\vp\,$ on a 
Riemannian manifold, provided that $\,v\,$ and $\,Q\,$ are still given by 
(\ref{not}). The identity (\ref{rcv}.a) is well known, cf.\ \cite{calabi}.
\begin{lemma}\label{fotcl}Let there be given a K\"ah\-ler surface\/ 
$\,(M,g)\,$ with the K\"ah\-ler form\/ $\,\kfo$, a Kil\-ling potential\/ 
$\,\vp\,$ on\/ $\,(M,g)$, and\/ $\,C^\infty$ functions\/ 
$\,\fx,\fy:M\to\bbR\hh$. Then, in the notation of\/ {\rm(\ref{not})}, the 
$\,2$-form\/ $\,\eta=\fx\hs\kfo\hs+\fy\,\xi\wedge\hh d\vp\,$ is closed if 
and only if\/ $\,d\hh(\fx-Q\hh\fy)+(d_v\fy+\fy Y)\hs d\vp
+(d_u\fy)\hh\xi=0$.
\end{lemma}
\begin{proof}For $\,A\,$ corresponding to $\,\eta\,$ as in 
Lemma~\ref{cltwf}, $\,JA=\fy(\xi\otimes u+d\vp\otimes v)-\fx$, where 
$\,\fx\,$ stands for $\,\fx\,$ times $\,\mathrm{Id}\hh$. Also, 
$\,\mathrm{div}\hs(\xi\otimes u)=\nabla_{\hskip-1.2ptu}\hh\xi\,$ and 
$\,\mathrm{div}\hs(d\vp\otimes v)=Yd\vp+\nabla_{\hskip-1.2ptv}\hh d\vp$, 
so that 
$\,\mathrm{div}\hs(\xi\otimes u+\hs d\vp\otimes v)=Yd\vp$. (Note that 
$\,\nabla_{\hskip-1.2ptu}\hh\xi=g(\nabla_{\hskip-1.2ptu}u,\,\cdot\,)$, 
which is the opposite of $\,\nabla_{\hskip-1.2ptv}\hh d\vp
=g(\nabla_{\hskip-1.2ptv}v,\,\cdot\,)$, as 
$\,\nabla_{\hskip-1.2ptu}u=\nabla_{\hskip-1.2ptu}(Jv)
=J\nabla_{\hskip-1.2ptu}v=J\nabla_{\hskip-1.2ptv}u
=\nabla_{\hskip-1.2ptv}(Ju)=-\nabla_{\hskip-1.2ptv}v$.) Thus, 
$\,\mathrm{div}\hs JA=\fy Yd\vp+(d_u\fy)\hh\xi+(d_v\fy)\hs d\vp-\df$. 
As $\,\lj\hs\omega,\eta\rg=2\fx-Q\hh\fy\,$ by (\ref{oab}), Lemma~\ref{cltwf} 
yields our claim.
\qed
\end{proof}
Given a nonconstant Kil\-ling potential $\,\vp\,$ on a compact K\"ah\-ler 
manifold $\,(M,g)\,$ and a $\,C^\infty$ function $\,\psi:M\to\bbR$, one may 
refer to $\,\psi\,$ as {\it a\/ $\,C^\infty$ function of\/} $\,\vp\,$ if 
$\,\psi=G(\vp)\,$ for some $\,C^\infty$ function 
$\,G:[\hh\vp_{\mathrm{min}},\vp_{\mathrm{max}}]\to\bbR\hh$. Note that
\begin{equation}\label{fcn}
\psi\hs\,\mathrm{\ is\ a\ }\,C^\infty\mathrm{\ function\ of\ }\,\hs\vp\hs\,
\mathrm{\ if\ and\ only\ if\ }\,\,d\psi\wedge d\vp=0\hh.
\end{equation}
In fact, let $\,M'\nh\subset M\,$ be the open set on which $\,d\vp\ne0$. It is 
well-known that $\,M'$ is connected and dense in $\,M$, and that Kil\-ling 
potentials are \hbox{Morse\hs-}\hskip0ptBott functions (cf.\ 
\cite[Remark~2.3(ii) and Example 11.1]{derdzinski-maschler-06}). The relation 
$\,d\psi\wedge d\vp=0\,$ clearly means that $\,\psi\,$ restricted to $\,M'$ 
is, locally, a $\,C^\infty$ function of $\,\vp$. Consequently, the word 
`locally' can be dropped, since the connectivity lemma, mentioned at the end 
of Section~\ref{pr}, now implies connectedness of the $\,\vp$-pre\-im\-ages of 
all real numbers. Also, due to the \hbox{Morse\hs-}\hskip0ptBott property of 
$\,\vp$, its critical manifolds are compact and isolated from one another, so 
that their number is finite, and, as $\,\vp\,$ is constant on each of them, 
$\,\vp\,$ has a finite set $\,\Gamma$ of critical values. Next, we show that 
the function 
$\,G:[\hh\vp_{\mathrm{min}},\vp_{\mathrm{max}}]\smallsetminus\Gamma\to\bbR\,$ 
with $\,\psi=G(\vp)\,$ on $\,M\smallsetminus\vp^{-1}(\Gamma)\,$ has a 
$\,C^\infty$ extension to $\,[\hh\vp_{\mathrm{min}},\vp_{\mathrm{max}}]$. To 
this end, we fix $\,\vp_*\in\Gamma\,$ and a point $\,x\in M\,$ 
such that $\,\vp(x)=\vp_*$ and $\,d\vp_x\nh=0$. The null\-space 
of the Hess\-i\-an of $\,\vp\,$ at $\,x\,$ coincides with the tangent space at 
$\,x\,$ of the critical manifold of $\,\vp\,$ containing $\,x\,$ (cf.\ 
\cite[Remark~2.3(iii-d)]{derdzinski-maschler-06}). One may thus choose 
$\,\delta,\ve\in(0,\infty)\,$ and a $\,C^\infty$ curve 
$\,(-\hh\ve,\ve)\ni t\mapsto x(t)\,$ in $\,M\,$ with $\,x(0)=x$, for which the 
assumptions, and hence the conclusion, of Lemma~\ref{cinff} are satisfied if 
one lets the symbols $\,\vp\,$ and $\,\psi\,$ stand for the functions 
$\,\vp(x(t))\,$ and $\,\psi(x(t))\,$ of the variable $\,t$.

\section{Special bi\-con\-for\-mal changes}\label{sb}
\setcounter{equation}{0}
Two Riemannian metrics $\,g,\hatg\,$ on a manifold $\,M\,$ are sometimes 
referred to as {\it bi\-con\-for\-mal\/} \cite{el-mansouri,ganchev-mihova} if 
there exist vector sub\-bundles $\,\mathcal{V}\,$ and $\,\mathcal{H}\,$ of 
$\,\tm\,$ with $\,\tm\nh=\mathcal{V}\oplus\mathcal{H}\,$ such that, for some 
positive $\,C^\infty$ functions $\,\fx,\chi:M\to\bbR$,
\begin{equation}\label{onh}
\hatg=\fx\hs g\hskip6pt\mathrm{on}\hskip4.5pt\mathcal{H}\hs,\hskip18pt
\hatg=\chi\hs g\hskip6pt\mathrm{on}\hskip4.5pt\mathcal{V},\hskip18pt
g(\mathcal{H},\mathcal{V})=\hatg(\mathcal{H},\mathcal{V})=\{0\}\hs.
\end{equation}
This kind of bi\-con\-for\-mal\-i\-ty is of little interest in the case of two 
K\"ah\-ler metrics on a given complex surface $\,M\nh$, since, locally, in a 
dense open subset of $\,M\nh$, (\ref{onh}) always holds, due to the existence 
of eigen\-space bundles of $\,\hatg\,$ relative to $\,g$.

The special bi\-con\-for\-mal changes defined in the Introduction represent a 
particular case of the situation described above. Namely, relation 
(\ref{bcm}.i), in the open set $\,M'\nh\subset M\,$ where $\,d\vp\ne0$, 
amounts to (\ref{onh}) with $\,\chi=\fx-Q\hh\fy$, 
$\,\mathcal{V}=\mathrm{Span}_{\bbR}(v,u)\,$ and 
$\,\mathcal{H}=\mathcal{V}^\perp$ (notation of (\ref{not})), so that 
$\,\tm'\nh=\mathcal{V}\oplus\mathcal{H}$, while
\begin{equation}\label{eig}
\fx\,\mathrm{\ and\ }\,\fx-Q\hh\fy\,\mathrm{\ are\ the\ eigenvalue\ functions\ 
of\ }\,\hatg\,\mathrm{\ relative\ to\ }\,g\hs.
\end{equation}

Given a K\"ah\-ler manifold $\,(M,g)\,$ with a nonconstant Kil\-ling potential 
$\,\vp\,$ and $\,C^\infty$ functions $\,\fx,\fy:M\to\bbR$, let a 
twice-co\-var\-i\-ant symmetric tensor field $\,\hatg\,$ on $\,M\,$ be 
Her\-mit\-i\-an relative to the underlying complex structure $\,J$. Thus, 
$\,\hato=\hatg(J\,\cdot\,,\,\cdot\,)\,$ is a $\,2$-form. Then (\ref{bcm}.i) 
holds, for $\,\xi=g(J(\navp),\,\cdot\,)$, if and only if 
\begin{equation}\label{hoe}
\hato\,\,=\,\,\fx\hs\kfo\hs\,\,+\,\,\fy\,\xi\wedge\hh d\vp\hs.
\end{equation}
A nontrivial special bi\-con\-for\-mal change (\ref{bcm}) of a 
met\-\hbox{ric-}\hskip.7ptpo\-ten\-tial pair $\,(g,\nnh\vp)$, if it exists, is never unique. 
Namely, it gives rise to a three-pa\-ram\-e\-ter family of such changes, 
leading to the met\-\hbox{ric-}\hskip.7ptpo\-ten\-tial pairs $\,(p\hatg+qg,p\havp+q\vp+s)$, 
with any constants $\,p,q,s\,$ such that $\,p\hatg+qg\,$ is positive 
definite (for instance, $\,p,q\,$ may both be positive). In fact, 
(\ref{bcm}) holds if one replaces $\,\hatg,\fx\nh,\fy\,$ and $\,\havp\,$ 
by $\,g\hh'\nh=p\hatg+qg,p\fx\nh+q,p\hs\fy\,$ and 
$\,\vp\hh'\nh=p\havp+q\vp+s$. (Specifically, (\ref{bcm}.ii) for $\,\vp\hh'$ is 
immediate as $\,\imath_vg\hh'\nh=p\imath_v\hatg+q\imath_vg$, for $\,v=\navp$.) 
In addition, $\,g\hh'$ is a K\"ah\-ler metric, since the $\,2$-form 
$\,\omega\hs'\nh=g\hh'(J\,\cdot\,,\,\cdot\,)\,$ equals 
$\,p\hh\hato\hs+q\hh\kfo$, and so $\,d\hs\omega\hh'\nh=0$.

The existence of a special bi\-con\-for\-mal change (\ref{bcm}) for a pair 
$\,(g,\nnh\vp)$, with prescribed $\,\fx\,$ and $\,\fy$, obviously amounts to 
requiring $\,\hatg\,$ given by (\ref{bcm}.i) to be a K\"ah\-ler metric such 
that $\,v=\navp\,$ is the $\,\hatg$-gra\-di\-ent of some $\,C^\infty$ function 
$\,\havp$. The following lemma describes a condition equivalent to this in the 
case of K\"ah\-ler surfaces; a similar result, valid in all complex 
dimensions, was obtained by Gan\-chev and Mi\-ho\-va 
\cite[the text following Definition 4.1]{ganchev-mihova}.
\begin{lemma}\label{spbic}Given a met\-\hbox{ric-}\hskip.7ptpo\-ten\-tial pair\/ $\,(g,\nnh\vp)\,$ 
on a compact complex surface\/ $\,M\,$ and\/ $\,C^\infty$ functions\/ 
$\,\fx,\fy,\havp:M\to\bbR$, one has\/ {\rm(\ref{bcm})} for a 
met\-\hbox{ric-}\hskip.7ptpo\-ten\-tial pair of the form\/ $\,(\hatg,\havp)\,$ on\/ $\,M\,$ 
if and only if, in the notation of\/ {\rm(\ref{not})},
\begin{enumerate}
  \def\theenumi{{\rm\roman{enumi}}}
\item $\havp=P(\vp)\,$ for some\/ $\,C^\infty$ function\/ 
$\,P:[\hh\vp_{\mathrm{min}},\vp_{\mathrm{max}}]\to\bbR$,
\item $\fx-\hh Q\hh\fy=\fc(\vp)$, with\/ $\,\fc=\hh dP/d\vp$,
\item $d_u\fy=0\,$ and\/ $\,d_v\fy+\fy Y\nh=-\hs\fc'(\vp)$, where\/ 
$\,\fc'=\hh d\fc/d\vp$,
\item $\fx>\hh\mathrm{max}\,(Q\hh\fy,0)$.
\end{enumerate}
Sufficiency of\/ {\rm(i)} -- {\rm(iv)} remains true without the compactness 
hypothesis.
\end{lemma}
\begin{proof}Necessity: first, (\ref{bcm}) implies (i). In fact, 
$\,d\havp=\imath_v\hatg=(\fx-Q\hh\fy)\hs d\vp\,$ in view of (\ref{bcm}), so 
that (\ref{fcn}) gives (i) and (ii). Next, by Lemma~\ref{fotcl}, (\ref{hoe}) 
and (ii), $\,[\hs d_v\fy+\fy Y\nh+\fc'(\vp)]\hs d\vp+(d_u\fy)\hh\xi=0$, and 
(iii) follows since $\,\xi\,$ is orthogonal to $\,d\vp$. Finally, (iv) amounts 
to positive definiteness of $\,\hatg$, cf.\ (\ref{eig}).

Sufficiency: conditions (ii) -- (iv) combined with Lemma~\ref{fotcl} and 
(\ref{hoe}) show that $\,\hatg\,$ defined by (\ref{bcm}.i) is a K\"ah\-ler 
metric. Also, in view of (\ref{bcm}.i) and (\ref{not}), 
$\,\imath_v\hatg=(\fx-Q\hh\fy)\hs d\vp$, which, according to (i) and (ii), 
equals $\,d\havp$. This proves (\ref{bcm}.ii).
\qed
\end{proof}
\begin{remark}\label{exloc}For any met\-\hbox{ric-}\hskip.7ptpo\-ten\-tial pair 
$\,(g,\nnh\vp)\,$ on a K\"ah\-ler surface, nontrivial special bi\-con\-for\-mal 
changes of $\,(g,\nnh\vp)\,$ exist locally, at points where $\,d\vp\ne0$, 
and the $\,C^\infty$ function $\,P\,$ of the variable $\,\vp$, such 
that $\,\havp=P(\vp)\,$ for the resulting pair $\,(\hatg,\havp)$, may be 
prescribed arbitrarily, as long as $\,dP/d\vp>0$. (Cf.\ 
Lemma~\ref{spbic}(i)-(ii) and (\ref{eig}).) This is clear from the final 
clause of Lemma~\ref{spbic}, since conditions (iii) and (iv) in 
Lemma~\ref{spbic} can be realized by solving an ordinary differential 
equation, with suitably chosen initial data, along each integral curve of 
$\,v$.
\end{remark}

\section{One general construction}\label{og}
\setcounter{equation}{0}
The following theorem provides a method of constructing examples of nontrivial 
special bi\-con\-for\-mal changes in complex dimension $\,2$. In the next four 
sections this method will be applied to four specific classes of K\"ah\-ler 
surface metrics.
\begin{theorem}\label{dftef}Given a nonconstant Kil\-ling potential\/ 
$\,\vp\,$ on a compact K\"ah\-ler surface\/ $\,(M,g)$, the following two 
conditions are equivalent\/{\rm:}
\begin{enumerate}
  \def\theenumi{{\rm\alph{enumi}}}
\item $(g,\nnh\vp)\,$ admits a nontrivial special bi\-con\-for\-mal change 
as in\/ {\rm(\ref{bcm})}, with\/ $\,\hs\fy\hs$ which is a\/ $\,C^\infty$ 
function of\/ $\,\vp$,
\item $\Delta\hh[S(\vp)]\hs=\hs-\hs\fc'(\vp)\,$ for some nonconstant\/ 
$\,C^\infty$ functions\/ 
$\,S,\fc:[\hh\vp_{\mathrm{min}},\vp_{\mathrm{max}}]\to\bbR\,$ and\/ 
$\,\fc'=\hh d\fc/d\vp$.
\end{enumerate}
Then, up to additive constants, $\,\fc\,$ in\hs\/ {\rm(b)} coincides with\/ 
$\,\fc\,$ appearing in Lemma~{\rm\ref{spbic}}, while\/ $\,\fy\,$ in\/ {\rm(a)} 
and\/ $\,S\,$ in\/ {\rm(b)} are related by\/ $\,\fy=dS/d\vp$.
\end{theorem}
\begin{proof}Assuming (a) and using Lemma~\ref{spbic}(iii), one obtains (b) 
for any $\,S\,$ with $\,dS/d\vp=\fy$. Conversely, (b) easily implies 
condition (iii) in Lemma~\ref{spbic}(iii) for $\,\fy=dS/d\vp$. Adding a 
suitable constant to $\,\fc\nh$, one also gets (iv) in Lemma~\ref{spbic} for 
$\,P\nh,\fx\,$ and $\,\havp\,$ chosen so as to satisfy (ii) and (i) in 
Lemma~\ref{spbic}.
\qed
\end{proof}

\section{K\"ah\-ler-Ein\-stein surfaces}\label{ke}
\setcounter{equation}{0}
On any compact K\"ah\-ler-Ein\-stein manifold $\,(M,g)\,$ such that the 
constant $\,\lambda\,$ with $\,\mathrm{Ric}\hh=\lambda\nh g\,$ is positive and 
$\,M\,$ admits a nontrivial hol\-o\-mor\-phic vector field, there exists a 
nonconstant Kil\-ling potential. In fact, by Matsu\-shi\-ma's theorem 
\cite{matsushima}, $\,\mathfrak{h}=\mathfrak{g}\oplus\nnh J\hn\mathfrak{g}\,$ 
for the spaces $\,\mathfrak{h}\,$ and $\,\mathfrak{g}\,$ of all 
real-hol\-o\-mor\-phic vector fields and, respectively, all 
real-hol\-o\-mor\-phic gradients, where $\,J\hn\mathfrak{g}\,$ consists of all 
Kil\-ling fields on $\,(M,g)$.

Using Theorem~\ref{dftef} one sees that {\it a nontrivial special 
bi\-con\-for\-mal change of $\,(g,\nnh\vp)\,$ exists whenever\/ $\,g\,$ is a 
K\"ah\-ler-Ein\-stein metric with positive Ein\-stein constant\/ 
$\,\lambda\,$ on a compact complex surface\/ $\,M\,$ and\/ 
$\,\vp\,$ is a nonconstant Kil\-ling potential on\/} $\,(M,g)$. Namely, by 
(\ref{rcv}.a), $\,\Delta\vp=a-2\lambda\nh\vp\,$ for some $\,a\in\bbR\hh$. 
Thus, condition (b) in Theorem~\ref{dftef} holds for $\,S(\vp)=\vp\,$ and 
$\,\fc(\vp)=\lambda\nh\vp^2\nh-a\vp$.

\section{K\"ah\-ler-Ric\-ci sol\-i\-tons}\label{kr}
\setcounter{equation}{0}
A {\it Ric\-ci soliton\/} \cite{hamilton} is a Riemannian manifold 
$\,(M,g)\,$ with the property that, for some constant $\,\lambda$, the tensor 
field $\,\lambda\nh gw-\hh\mathrm{Ric}\,$ is the Lie derivative of $\,g\,$ in 
the direction of some vector field. Perelman \cite[Remark~3.2]{perelman} 
proved that, if $\,M\,$ is compact, such a vector field must be the sum of a 
Kil\-ling field and a gradient, or, equivalently, there exists a $\,C^\infty$ 
function $\,\vp:M\to\bbR\,$ with
\begin{equation}\label{ndt}
\nabla d\vp\,+\,\hh\mathrm{Ric}\hh\,=\,\lambda\nh g\hskip9pt\mathrm{for\ a\ 
constant\ }\,\lambda\hs. 
\end{equation}
If a $\,C^\infty$ function $\,\vp:M\to\bbR\,$ satisfies (\ref{ndt}), then 
\cite[p.\ 201]{chow}
\begin{equation}\label{cte}
c\,=\,\Delta\vp\,-\,g(\navp,\navp)\,+\,2\lambda\nh\vp
\end{equation}
is a constant. In fact, adopting the notation of (\ref{not}) except for the 
formulae involving $\,J$, and applying to both sides of (\ref{ndt}) either 
$-2\,\mathrm{div}$, or $\,\mathrm{tr}_g$ followed by $\,d$, or, finally, 
$\,2\hs\imath_v$, one obtains 
$\,-2\hs d\hskip.2ptY\nh-2\hskip1.5pt\mathrm{Ric}\hh(v,\,\cdot\,)-d\hs\se=0$, 
$\,d\hskip.2ptY\nh+d\hs\se=0$, and, by (\ref{rcv}.b), 
$\,dQ+2\hskip1.5pt\mathrm{Ric}\hh(v,\,\cdot\,)=2\lambda\hs d\vp$. (Here 
$\,\se\,$ denotes the scalar curvature, while 
$\,\mathrm{div}\nabla d\vp=d\hskip.2ptY\nh+\mathrm{Ric}\hh(v,\,\cdot\,)\,$ by 
the Boch\-ner identity, which has the coordinate form 
$\,v^{\hs k}{}_{,\hs jk}\nh=v^{\hs k}{}_{,\hs kj}+R_{jk}v^{\hs k}\nnh$, and 
$\,\mathrm{div}\,\mathrm{Ric}\hh=d\hs\se\hs/2\,$ in 
view of the Bian\-chi identity for the Ric\-ci tensor.) Adding these three 
equalities produces the relation 
$\,d\hh[\Delta\vp-g(\navp,\navp)+2\lambda\nh\vp]=0$.

By a {\it K\"ah\-ler-Ric\-ci sol\-i\-ton\/} one means a Ric\-ci sol\-i\-ton 
which is at the same time a K\"ah\-ler manifold \cite{tian-zhu,wang-zhu}. A 
function $\,\vp\,$ with (\ref{ndt}) then is a Kil\-ling potential. (Since 
$\,g\,$ and $\,\mathrm{Ric}\,$ are Her\-mit\-i\-an, so must be 
$\,\nabla d\vp\,$ as well.) Also, (\ref{cte}) can be rewritten as 
$\,\Delta\emvp\nh=(2\lambda\nh\vp-c)\hh\emvp\nnh$. Thus, by 
Theorem~\ref{dftef}, {\it for every non-Ein\-stein compact K\"ah\-ler-Ric\-ci 
sol\-i\-ton\/ $\,(M,g)\,$ of complex dimension\/ $\,2$, the pair\/ 
$\,(g,\nnh\vp)$, where\/ $\,\vp\,$ is a function satisfying\/ 
{\rm(\ref{ndt})}, admits a nontrivial special bi\-con\-for\-mal change}. 
Specifically, condition (b) in Theorem~\ref{dftef} then holds for 
$\,S(\vp)=\emvp$ and $\,\fc(\vp)=[2\lambda\hh(\vp+1)-c\hh]\hh\emvp\nnh$.

\section{Con\-for\-mal\-ly-Ein\-stein K\"ah\-ler surfaces}\label{ce}
\setcounter{equation}{0}
Let $\,(M,g)\,$ be a con\-for\-mal\-ly-Ein\-stein, non-Ein\-stein compact 
K\"ah\-ler surface. The scalar curvature $\,\se\,$ of $\,g\,$ then is a 
nonconstant Kil\-ling potential, and so $\,g\,$ is an extremal K\"ah\-ler 
metric \cite{calabi}, while $\,\se>0\,$ everywhere and 
$\,\hs\se\hh^3\nh+6\hs\se\hs Y\nh-12\hs Q=12\hh\cj\,$ for some constant 
$\,\cj>0$. See \cite[Proposition 4 on p.\ 419 and Theorem 2 on p.\ 
428]{derdzinski-83}, \cite[Lemma 3 on p. 169]{lebrun}.

Con\-for\-mal\-ly-Ein\-stein, non-Ein\-stein K\"ah\-ler metrics are known to 
exist on both the one-point and two-point blow-ups of $\,\bbCP^2\nnh$. The 
former, found by Ca\-la\-bi \cite{calabi}, is conformal to the Page metric 
\cite{page}, for reasons given in \cite[the top of p.\ 430]{derdzinski-83}; 
the existence of the latter is a result of Chen, LeBrun and Weber 
\cite{chen-lebrun-weber}.

Theorem~\ref{dftef} implies that {\it for every 
con\-for\-mal\-ly-Ein\-stein, non-Ein\-stein compact K\"ah\-ler surface\/ 
$\,(M,g)$, the pair\/ $\,(g,\nnh\vp)$, with $\,\vp=\se\hh$, admits a nontrivial 
special bi\-con\-for\-mal change}. In fact, the equality 
$\,\hs\se\hh^3\nh+6\hs\se\hs Y\nh-12\hs Q=12\hh\cj\,$ yields condition (b) in 
Theorem~\ref{dftef} for $\,\vp=\se\,$ and $\,S(\vp)=-\hh\vp^{-1}\nnh$, with 
$\,\fc(\vp)=\cj\hskip1pt\vp^{-2}\nh+\hh\vp/6$. The existence of such a 
bi\-con\-for\-mal change in this case was first discovered by LeBrun 
\cite[p.\ 171, the end of the proof of Proposition 2]{lebrun}, who 
proved that, with $\,\rho\,$ and $\,\kfo\,$ standing for the Ric\-ci and 
K\"ah\-ler forms of $\,g$,
\begin{equation}\label{pos}
\rho\,\,+\,\,2\hs i\hskip1pt\partial\overline{\partial}\hskip1pt\log\hs\se\hs\,
\,\,=\,\,[(Q+\cj)\hs\se\hh^{-2}+\,\se/6\hh]\,\omega\,\,
+\,\,\hs\se^{-2}\hs\xi\wedge\hh d\hs\se
\end{equation}
(notation of (\ref{not}). Equality (\ref{pos}) easily implies (\ref{bcm}.i) 
with a new K\"ah\-ler metric $\,\hat g$. Namely, the right-hand side of 
(\ref{pos}) coincides with $\,\hato\,$ in (\ref{hoe}), for suitable $\,\fx\,$ 
and $\,\fy$, while the Her\-mit\-i\-an $\,2$-ten\-sor field $\,\hatg\,$ 
characterized by $\,\hato=\hatg(J\,\cdot\,,\,\cdot\,)\,$ is positive definite, 
cf.\ Lemma~\ref{spbic}(iv); at the same time, the left-hand side of 
(\ref{pos}) is a closed $\,2$-form.

\section{Special K\"ah\-ler-Ric\-ci potentials}\label{sk}
\setcounter{equation}{0}
A {\it special K\"ah\-ler-Ric\-ci potential\/} \cite{derdzinski-maschler-03} 
on a compact K\"ah\-ler surface $\,(M,g)\,$ is any nonconstant Kil\-ling 
potential $\,\vp:M\to\bbR\,$ such that both $\,Q=g(\navp,\navp)\,$ and 
$\,Y\nnh=\Delta\vp\,$ are $\,C^\infty$ functions of $\,\vp$. This definition, 
although different from the one given in 
\cite[\S~7]{derdzinski-maschler-03} for all complex dimensions $\,m\ge2$, 
is equivalent to it when $\,m=2$, as a consequence of (\ref{rcv}) and 
(\ref{fcn}). Special K\"ah\-ler-Ric\-ci potentials on compact K\"ah\-ler 
manifolds $\,(M,g)\,$ of any given complex dimension $\,m\ge2\,$ were 
classified in \cite{derdzinski-maschler-06}. They turn out to be 
bi\-hol\-o\-mor\-phic to $\,\bbCP^m$ or to hol\-o\-mor\-phic $\,\bbCP^1$ 
bundles over complex manifolds admitting K\"ah\-ler-Ein\-stein metrics. If 
$\,m\ge3$, there are other natural conditions which imply the existence of a 
K\"ah\-ler-Ric\-ci potential \cite[Corollary 9.3]{derdzinski-maschler-03}, 
\cite[Theorem 6.4]{jelonek}.

If $\,(M,g)\,$ is a compact K\"ah\-ler surface with a special 
K\"ah\-ler-Ric\-ci potential $\,\vp$, {\it any\/} nonconstant $\,C^\infty$ 
function $\,S:[\hh\vp_{\mathrm{min}},\vp_{\mathrm{max}}]\to\bbR\,$ satisfies 
condition (b) in Theorem~\ref{dftef}. Consequently, by Theorem~\ref{dftef}, 
{\it the pair\/ $\,(g,\nnh\vp)\,$ then admits a nontrivial special 
bi\-con\-for\-mal change\/} (\ref{bcm}), {\it in which\/ $\,\fy\,$ may be any 
prescribed\/ $\,C^\infty$ function of\/ $\,\vp\,$ other than the zero 
function.}

\section{Geodesic gradients}\label{gg}
\setcounter{equation}{0}
We say that a nonconstant Kil\-ling potential $\,\vp\,$ on a compact 
K\"ah\-ler manifold $\,(M,g)\,$ has a {\it geodesic gradient\/} if all the 
integral curves of $\,\navp\,$ are re\-pa\-ram\-e\-trized geodesics. By 
(\ref{fcn}), this amounts to requiring that $\,Q=g(\navp,\navp)\,$ be a 
$\,C^\infty$ function of $\,\vp$, since (\ref{rcv}.b) gives 
$\,2\hh\nabla_{\hskip-1.2ptv}v=\nabla Q\,$ (notation of (\ref{not})).

Thus, every special K\"ah\-ler-Ric\-ci potential (Section~\ref{sk}) has a 
geodesic gradient. Further examples, which are not special K\"ah\-ler-Ric\-ci 
potentials, are described in the Appendix.
\begin{theorem}\label{geodg}For every nonconstant Kil\-ling potential\/ 
$\,\vp\,$ with a geodesic gradient on a compact K\"ah\-ler surface\/ 
$\,(M,g)$, other than a special K\"ah\-ler-Ric\-ci potential, there exists a 
nontrivial special bi\-con\-for\-mal change\/ {\rm(\ref{bcm})} of the pair\/ 
$\,(g,\nnh\vp)$, for which\/ $\,\fc\nh$, defined in Lemma~{\rm\ref{spbic}}, 
can be, up to an additive constant, any prescribed nonconstant\/ $\,C^\infty$ 
function of the variable\/ 
$\,\vp\in[\hh\vp_{\mathrm{min}},\vp_{\mathrm{max}}]\,$ such that\/ 
$\,\fc'=\hh d\fc/d\vp\,$ is\/ $\,L^2$-or\-thog\-o\-nal to linear functions 
of\/ $\,\vp$.
\end{theorem}
For a proof, see the final paragraph of the Appendix.

\section{Bi\-con\-for\-mal changes defined on an open 
sub\-man\-i\-fold}\label{os}
\setcounter{equation}{0}
The last five sections described examples of nontrivial special 
bi\-con\-for\-mal changes that naturally arise in certain classes of compact 
K\"ah\-ler surfaces. As we will see below and in Section~\ref{ac}, there are 
also circumstances in which, for a given 
met\-\hbox{ric-}\hskip.7ptpo\-ten\-tial pair $\,(g,\nnh\vp)\,$ on a compact 
complex surface $\,M\nh$, one naturally obtains a nontrivial special 
bi\-con\-for\-mal change (\ref{bcm}) of $\,(g,\nnh\vp)\,$ restricted to the 
dense open sub\-man\-i\-fold $\,M'$ characterized by the condition 
$\,d\vp\ne0$, while the functions $\,\fx,Q\hh\fy,\havp\,$ in (\ref{bcm}), cf.\ 
(\ref{not}), and the metric $\,\hatg$, all have $\,C^\infty$ extensions to 
$\,M\nh$. The only difference between this case and the standard one (defined 
in the Introduction) is that $\,\fy$, unlike $\,Q\hh\fy$, may now fail to have 
a $\,C^\infty$ extension to $\,M\nh$.

To provide an example of such a situation, we let $\,M\,$ stand either for 
$\,\bbCP^2$ or for the one-point blow-up of $\,\bbCP^2\nnh$, so that $\,M\,$ 
is a simply connected compact complex surface with an effective action of 
$\,\mathrm{U}(2)\,$ by bi\-hol\-o\-mor\-phisms. For any 
$\,\mathrm{U}(2)$-in\-var\-i\-ant K\"ah\-ler metric $\,g\,$ on $\,M\nh$, a 
fixed vector field $\,u\,$ generating the action of the center 
$\,\mathrm{U}(1)\subset\hs\mathrm{U}(2)\,$ is a 
$\,\mathrm{U}(2)$-in\-var\-i\-ant $\,g$-Kil\-ling field. Thus, 
$\,u=J(\navp)\,$ for some nonconstant Kil\-ling potential $\,\vp\,$ on 
$\,(M,g)$. (Cf.\ \cite[Lemma~5.3]{derdzinski-maschler-03}.) One may say that 
such $\,\vp\,$ is {\it associated with\/} $\,g$. As the principal orbits of 
the $\,\mathrm{U}(2)\,$ action are three-di\-men\-sion\-al, every\/ 
$\,\mathrm{U}(2)$-in\-var\-i\-ant $\,C^\infty$ function $\,M\to\bbR\,$ is, by 
(\ref{fcn}), a $\,C^\infty$ function of $\,\vp$. Applied to the functions 
$\,Q=g(\navp,\navp)\,$ and $\,Y\nnh=\Delta\vp$, this shows that $\,\vp\,$ then 
is a special K\"ah\-ler-Ric\-ci potential on $\,(M,g)$, cf.\ Section~\ref{sk}.
\begin{lemma}\label{utinv}For\/ $\,M\,$ and\/ $\,M'\nh\subset M\,$ as above, 
any\/ $\,\mathrm{U}(2)$-in\-var\-i\-ant K\"ah\-ler metrics\/ $\,g,\hatg\,$ 
on\/ $\,M\nh$, and nonconstant Kil\-ling potentials\/ $\,\vp,\havp\,$ 
associated with them, the pair\/ $\,(\hatg,\havp)\,$ restricted to\/ $\,M'$ 
arises from\/ $\,(g,\nnh\vp)\,$ by a bi\-con\-for\-mal change\/ 
{\rm(\ref{bcm})}. If, in addition, $\,M\,$ is the one-point blow-up of\/ 
$\,\bbCP^2\nnh$, then the functions\/ $\,\fx\,$ and\/ $\,Q\hh\fy\,$ in\/ 
{\rm(\ref{bcm})}, as well as the corresponding eigen\-space bundles\/ 
$\,\mathcal{V}\,$ and\/ $\,\mathcal{H}$, introduced in Section~{\rm\ref{sb}}, 
all have\/ $\,C^\infty$ extensions to\/ $\,M\,$ such that\/ 
$\,\fx>\hh\mathrm{max}\,(Q\hh\fy,0)\,$ on\/ $\,M\nh$, cf.\ {\rm(\ref{not})}.
\end{lemma} 
\begin{proof}Let us denote by $\,\mathcal{V}\,$ the com\-plex-line sub\-bundle 
of $\,\tm'\nnh$, spanned by $\,v=\navp\,$ (that is, by $\,u=Jv$). The 
$\,g$-or\-thog\-o\-nal complement $\,\mathcal{H}=\mathcal{V}^\perp$ of 
$\,\mathcal{V}\,$ in $\,\tm'\nnh$, relative to any 
$\,\mathrm{U}(2)$-in\-var\-i\-ant K\"ah\-ler metric $\,g$, does not depend on 
the choice of such a metric. In fact, being $\,g$-or\-thog\-o\-nal to 
$\,v=-Ju$, the sub\-bundle $\,\mathcal{H}\,$ is contained in the real 
three-di\-men\-sion\-al sub\-bundle tangent to the orbits of the 
$\,\mathrm{U}(2)\,$ action. (As $\,v\,$ is the $\,g$-gra\-di\-ent of the 
$\,\mathrm{U}(2)$-in\-var\-i\-ant function $\,\vp$, it is orthogonal to the 
orbits.) Independence of $\,\mathcal{H}\,$ from $\,g\,$ now follows since a 
real three-di\-men\-sion\-al subspace of a complex two-di\-men\-sion\-al 
vector space contains only one complex subspace of complex dimension $\,1$.

For any two $\,\mathrm{U}(2)$-in\-var\-i\-ant K\"ah\-ler metrics $\,g,\hatg$, 
one clearly has (\ref{onh}) with some positive $\,C^\infty$ functions 
$\,\fx,\chi:M'\nh\to\bbR\hh$. Defining  $\,\fy$, on $\,M'\nnh$, by 
$\,\chi=\fx-Q\hh\fy$, we now obtain (\ref{bcm}.i) on $\,M'\nnh$, while 
(\ref{bcm}.ii) is obvious as $\,\hanb\havp=-Ju=\nabla\vp$. Finally, if $\,M\,$ 
is the one-point blow-up of $\,\bbCP^2\nnh$, then $\,\mathcal{V}\,$ (and hence 
$\,\mathcal{H}\,$ as well) has a $\,C^\infty$ extension to a sub\-bundle of 
$\,\tm\nh$, due to the fact that $\,M\,$ is a hol\-o\-mor\-phic $\,\bbCP^1$ 
bundle over $\,\bbCP^1$ and $\,\mathcal{V}\,$ is tangent to the fibres.
\qed
\end{proof}

\section{More on $\,\mathrm{U}(2)$-in\-var\-i\-ant K\"ah\-ler 
metrics}\label{ut}
\setcounter{equation}{0}
Let us consider the one-point blow-up $\,M\,$ of $\,\bbCP^2\nnh$, with the 
effective action of $\,\mathrm{U}(2)\,$ by bi\-hol\-o\-mor\-phisms. By {\it 
central au\-to\-mor\-phisms\/} of $\,M\,$ we mean transformations that belong 
the hol\-o\-mor\-phic action of $\,\bbC^*$ on $\,M\nh$, generated by the 
action of the center $\,\mathrm{U}(1)\subset\hs\mathrm{U}(2)$. They all 
commute with the action of $\,\mathrm{U}(2)$. Thus, the pull\-back of any 
$\,\mathrm{U}(2)$-in\-var\-i\-ant K\"ah\-ler metric on $\,M\,$ under any 
central au\-to\-mor\-phism of $\,M\,$ is again a 
$\,\mathrm{U}(2)$-in\-var\-i\-ant K\"ah\-ler metric on $\,M\nh$.

Suppose, in addition, that $\,g,\hatg\,$ are $\,\mathrm{U}(2)$-in\-var\-i\-ant 
K\"ah\-ler metrics on $\,M$, and let $\,\vp,\havp\,$ denote the nonconstant 
Kil\-ling potentials associated with them (Section~\ref{os}). According to 
the final clause Lemma~\ref{utinv} and (\ref{eig}), one has (\ref{onh}) with 
positive $\,C^\infty$ functions $\,\fx,\chi:M\to\bbR$, where 
$\,\chi=\fx-Q\hh\fy$. Both $\,\fx\,$ and $\,\chi\,$ are constant on either of 
the two exceptional orbits $\,\Sigma^\pm$ of the $\,\mathrm{U}(2)\,$ action, 
bi\-hol\-o\-mor\-phic to $\,\bbCP^1\nnh$. (In fact, $\,\fx\,$ and $\,\chi\,$ 
are $\,\mathrm{U}(2)$-in\-var\-i\-ant, since so are both metrics.) Let the 
constants $\,\fx^\pm$ and $\,\chi^\pm$ be the values of $\,\fx\,$ and 
$\,\chi\,$ on $\,\Sigma^\pm\nnh$. The positive real number 
\begin{equation}\label{dgg}
\mathrm{d}\hs(g,\nnh\hatg)\,=\,\chi^+\chi^-\nh/(\fx^+\fx^-)
\end{equation}
is an invariant which remains unchanged when one of the metrics $\,g,\hatg\,$ 
is replaced with its pull\-back under any central au\-to\-mor\-phism of 
$\,M\,$ (since the pair $\,(\chi^+\nnh,\chi^-)\,$ then is replaced by 
$\,(r\chi^+\nnh,r^{-1}\chi^-)\,$ for some $\,r\in(0,\infty)$). On the other 
hand, $\,\mathrm{d}\hs(g,\nnh\hatg)=1\,$ when $\,(\hatg,\havp)\,$ arises from 
$\,(g,\nnh\vp)\,$ by a special bi\-con\-for\-mal change: in fact, 
$\,\fx^\pm\nnh=\chi^\pm\nnh$, as $\,\chi=\fx-Q\hh\fy\,$ and $\,Q=0\,$ on 
$\,\Sigma^\pm\nnh$.

This shows that (b) implies (a) in the following proposition.
\begin{proposition}\label{utsbc}For any\/ $\,\mathrm{U}(2)$-in\-var\-i\-ant 
K\"ah\-ler metrics\/ $\,g,\hatg\,$ on the one-point blow-up of\/ 
$\,\bbCP^2\nnh$, and nonconstant Kil\-ling potentials\/ $\,\vp,\havp\,$ 
associated with them in the sense of Section~{\rm\ref{os}}, the following two 
conditions are equivalent\/{\rm:}
\begin{enumerate}
  \def\theenumi{{\rm\alph{enumi}}}
\item $\mathrm{d}\hs(g,\nnh\hatg)\hs=\hs1$,
\item $(\hatg,\havp)\,$ arises 
%from\/ $\,(g,\nnh\vp)\,$ by a special bi\-con\-for\-mal change.
by a special bi\-con\-for\-mal change from the pull\-back of\/ $\,(g,\nnh\vp)\,$ 
under some central au\-to\-mor\-phism of $\,M\nh$.
\end{enumerate}
\end{proposition} 
\begin{proof}It suffices to verify that (a) leads to (b). Because of how 
$\,\chi^\pm$ change under the action of a central au\-to\-mor\-phism (see 
above), one may use a pull\-back as in (b) to replace $\,\chi^+$ with the 
value $\,\fx^+\nnh$. As $\,\mathrm{d}\hs(g,\nnh\hatg)\hs=\hs1$, (\ref{dgg}) now 
gives $\,\chi^-\nh=\fx^-$ as well. Since $\,\chi=\fx-Q\hh\fy$, it follows that 
$\,Q\hh\fy=0\,$ on $\,\Sigma^\pm\nnh$.

In view of the final clause of Lemma~\ref{utinv}, the assertion will follow if 
one shows that $\,\fy\,$ (and not just $\,Q\hh\fy$) has a $\,C^\infty$ 
extension to $\,M\nh$. To this end, fix a point $\,x\in\Sigma^\pm$ and 
identify a neighborhood of $\,x\,$ in $\,M\,$ dif\-feo\-mor\-phic\-al\-ly with 
$\,\,U\nh\times D$, so that the flow of $\,u\,$ consists of the rotations 
$\,(y,z)\mapsto(y,zq)$, where $\,q\in\bbC\,$ and $\,|q|=1\,$ (notation of 
Remark~\ref{divis}, for $\,k=2$). According to Remark~\ref{divis}, both 
$\,Q\,$ and $\,Q\hh\fy\,$ is smoothly divisible by $\,|z|^2\nnh$, while 
$\,Q/|z|^2$ is positive on $\,\,U\nh\times\{0\}$, and so 
$\,\fy=[(Q\hh\fy)/|z|^2]/(Q/|z|^2)\,$ is smooth everywhere in 
$\,\,U\nh\times D$. That the Hess\-i\-an of $\,Q\,$ is nonzero everywhere in 
$\,\Sigma^\pm$ follows since the same is true of the Hessian of $\,d\vp$, cf.\ 
\cite[Remark 5.4]{derdzinski-maschler-03}. Namely, differentiating 
(\ref{rcv}.b) one sees that the former Hess\-i\-an equals twice the square of 
the latter, if both are identified with morphisms $\,\tm\to\tm\,$ as in the 
lines preceding (\ref{oab}).
\qed
\end{proof}
Let $\,g\,$ and $\,\hatg\,$ denote the two distinguished 
$\,\mathrm{U}(2)$-in\-var\-i\-ant K\"ah\-ler metrics on the one-point 
blow-up of $\,\bbCP^2\nnh$, mentioned at the end of the Introduction. 
According to Lemma~\ref{utinv}, the corresponding pairs $\,(g,\nnh\vp)\,$ and 
$\,(\hatg,\havp)\,$ arise from each other by the weaker version of a 
bi\-con\-for\-mal change, described at the beginning of Section~\ref{os}. The 
value of $\,\mathrm{d}\hs(g,\nnh\hatg)\,$ in this case is not known; if that value 
turns out to be  $\,1$, a stronger conclusion will be immediate from 
Proposition~\ref{utsbc}.

\section{Another construction}\label{ac}
\setcounter{equation}{0}
In contrast with Theorem~\ref{dftef}, the following result may lead to 
bi\-con\-for\-mal changes of a more general kind, introduced at the 
beginning of in Section~\ref{os}.
\begin{theorem}\label{opidd}Suppose that\/ $\,\vp\,$ is a nonconstant Kil\-ling 
potential on a compact K\"ah\-ler surface\/ $\,(M,g)\,$ and a K\"ah\-ler 
metric\/ $\,\hatg\,$ on\/ $\,M\,$ represents the same K\"ah\-ler cohomology 
class as\/ $\,g$. Using the notation of\/ {\rm(\ref{not})}, let us fix a\/ 
$\,C^\infty$ function\/ $\,\psi:M\to\bbR\,$ such that the K\"ah\-ler 
forms\/ $\,\kfo\,$ of\/ $\,g\,$ and\/ $\,\hato\,$ of\/ $\,\hatg\,$ are related 
by\/ $\,\hato=\kfo\hs+2i\hskip1pt\partial\overline{\partial}\hs\psi$, and 
denote by\/ $\,M'$ the open subset of\/ $\,M\,$ on which\/ $\,d\vp\ne0$.

If there exists a special bi\-con\-for\-mal change of\/ 
$\,(g,\nnh\vp)\,$ leading to a pair\/ $\,(\hatg,\havp)$, for some nonconstant 
Kil\-ling potential\/ $\,\havp\,$ on\/ $\,(M,\hatg)$, then
\begin{equation}\label{dvf}
d_v\psi\,\,\hs\mathrm{\ is\ a\ }\,C^\infty\mathrm{\ function\ of\ }\,\,\vp\hs\,
\mathrm{\ and\ }\hs\,d_u\psi=0\hs.
\end{equation}

Conversely, if\/ {\rm(\ref{dvf})} holds, then, for some nonconstant 
Kil\-ling potential\/ $\,\havp\,$ on\/ $\,(M,\hatg)$, the pair\/ 
$\,(\hatg,\havp)$, restricted to\/ $\,M'\nnh$, arises from\/ $\,(g,\nnh\vp)\,$ by 
a special bi\-con\-for\-mal change on\/ $\,M'\nnh$, and, on\/ 
$\,M'\nnh$, one has\/ {\rm(\ref{bcm})} with
\begin{equation}\label{fed}
\fx=\Delta\psi+1-\hh d(d_v\psi)/d\vp\hs,\hskip9pt
\fy=[\Delta\psi-2\hh d(d_v\psi)/d\vp]/Q\hs,\hskip9pt\havp=\vp+d_v\psi\hs.
\end{equation}
\end{theorem}
\begin{proof}Let some special bi\-con\-for\-mal change, applied to 
$\,(g,\nnh\vp)$, produce $\,(\hatg,\havp)$. Since $\,u\,$ is a Kil\-ling field for 
both $\,g\,$ and $\,\hatg\,$ (Section~\ref{pr}), the Lie derivatives 
$\,\Lie_u\kfo\,$ and $\,\Lie_u\hato\,$ both vanish, 
while $\,\Lie_u\,$ commutes with $\,\partial\overline{\partial}\,$ 
as $\,u\,$ is hol\-o\-mor\-phic. Thus, $\,\,d_u\psi\,$ lies in the kernel of 
$\,\partial\overline{\partial}\,$ and vanishes at points where $\,d\vp=0$, 
which is only possible if $\,\,d_u\psi=0\,$ identically. Now, by 
(\ref{bcm}.ii) and Lemma~\ref{hnteq}, $\,d\havp=d\vp+d\hs(d_v\psi)$. 
Consequently, Lemma~\ref{spbic}(i) yields (\ref{dvf}).

Conversely, let us assume (\ref{dvf}) and define $\,\fx,\fy,\havp\,$ by 
(\ref{fed}). Lemma~\ref{hnteq} gives 
$\,\hatg(v,\,\cdot\,)-g(v,\,\cdot\,)=d\havp-d\vp=d(d_v\psi)$, so that 
$\,\hatg(v,\,\cdot\,)\,$ equals the function 
$\,1+\hh d(d_v\psi)/d\vp=\fx-Q\hh\fy\,$ times 
$\,d\vp=g(v,\,\cdot\,)$. Hence $\,v\,$ is, at every point of $\,M'\nnh$, an 
eigenvector, for the eigenvalue $\,\fx-Q\hh\fy$, of $\,\hatg\,$ treated, with 
the aid of $\,g$, as a bundle morphism $\,\tm'\nh\to\tm'\nnh$. The 
Her\-mit\-i\-an $\,2$-ten\-sor field 
$\,\pi=(d\vp\otimes d\vp\hs+\hs\xi\otimes\xi)/Q\,$ 
(notation of (\ref{not})) is, obviously, the orthogonal projection onto the 
com\-plex-line sub\-bundle $\,\mathcal{V}\,$ of $\,\tm'\nnh$, spanned by $\,v$,
provided that one identifies $\,\pi$, as in the lines preceding (\ref{oab}), 
with a morphism $\,A:\tm'\nh\to\tm'\nnh$. Similarly, $\,g-\pi\,$ is the
orthogonal projection onto $\,\mathcal{H}=\mathcal{V}^\perp\nnh$. The other 
eigenvalue of $\,\hatg$, corresponding to eigenvectors in $\,\mathcal{H}$, is 
$\,\fx\nh$. (In fact, the sum of the two eigenvalues is 
$\,\mathrm{tr}_g/2$, which equals $\,2+\Delta\psi$, as one sees noting that, 
by (\ref{idd}), the relation 
$\,\hato=\kfo+2i\hskip1pt\partial\overline{\partial}\hs\psi\,$ amounts to 
$\,\hatg=g+\nabla d\psi+(\nabla d\psi)(J\,\cdot,J\,\cdot\,)$.) The spectral 
decomposition $\,\hatg=\fx(g-\pi)+(\fx-Q\hh\fy)\hs\pi\,$ now implies 
(\ref{bcm}.i), on $\,M'\nnh$, while Lemma~\ref{hnteq} yields (\ref{bcm}.ii), 
completing the proof.
\qed
\end{proof}

\section{The integral obstruction}\label{io}
\setcounter{equation}{0}
One can ask whether a nontrivial special bi\-con\-for\-mal changes 
exists for every met\-\hbox{ric-}\hskip.7ptpo\-ten\-tial pair $\,(g,\nnh\vp)\,$ on a compact 
complex surface $\,M\nh$. Here are two comments related to this existence 
question.

First, due to compactness of $\,M\nh$, for such a special bi\-con\-for\-mal 
change (\ref{bcm}), the functions $\,\fx\,$ and $\,\fy\,$ are uniquely 
determined by 
$\,\fc:[\hh\vp_{\mathrm{min}},\vp_{\mathrm{max}}]\to\bbR\,$ appearing in 
Lemma~\ref{spbic}. Namely, the zero 
function is the only $\,C^\infty$ solution $\,\phi:M\to\bbR\,$ to the 
homogeneous equation $\,d_v\phi+\phi Y=0$, associated with the equation 
imposed on $\,\fy\,$ in Lemma~\ref{spbic}(iii). In fact, the Kil\-ling 
potential $\,\vp\,$ has a finite number of critical manifolds (see the lines 
following (\ref{fcn})). At the same time, $\,Y\nnh=\Delta\vp\,$ is negative at 
points where $\,\vp=\vp_{\mathrm{max}}$, since the Hess\-i\-an 
$\,\nabla d\vp\,$ is negative sem\-i\-def\-i\-nite and nonzero at such points 
(cf.\ \cite[Remark 5.4]{derdzinski-maschler-03}). Thus, there exist 
$\,\delta,\ve\in(0,\infty)\,$ with the property that $\,d\vp\ne0\,$ and 
$\,Y\nh\le-\hs\delta\,$ everywhere in the open set $\,\,U\,$ on which 
$\,0<\vp_{\mathrm{max}}\nh-\vp<\ve$. Any integral curve 
$\,[\hs0,\infty)\ni t\mapsto x(t)\,$ of $\,v=\navp\,$ with 
$\,\vp(x(0))\in U\,$ lies entirely in $\,\,U\nh$. A solution $\,\phi\,$ to 
$\,d_v\phi+\phi Y=0$, if not identically zero along the integral curve, 
may be assumed positive everywhere on the curve, and then 
$\,d\hs[\hh\log\hh\phi(x(t))]/dt=-Y(x(t))\ge\delta>0$, so that 
$\,\phi(x(t))\to\infty\,$ as $\,t\to\infty$, contrary to compactness of 
$\,M\nh$. Thus, $\,\fc\,$ determines $\,\fy$, and hence $\,\fx\nnh$, 
cf.\ Lemma~\ref{spbic}(ii).

Secondly, let us fix a met\-\hbox{ric-}\hskip.7ptpo\-ten\-tial pair $\,(g,\nnh\vp)\,$ on a compact 
complex surface $\,M\nh$. In an attempt to find a nontrivial special 
bi\-con\-for\-mal change of $\,(g,\nnh\vp)$, one might begin by selecting a 
nonconstant $\,C^\infty$ function 
$\,\fc:[\hh\vp_{\mathrm{min}},\vp_{\mathrm{max}}]\to\bbR$, which would then 
become the function $\,\fc\,$ corresponding to such a bi\-con\-for\-mal change 
as in Lemma~\ref{spbic}(ii). (It is nonconstant as we want the change to be 
nontrivial, cf.\ the preceding paragraph.) Using the notation of (\ref{not}), 
we consider an arbitrary maximal integral curve $\,\bbR\ni t\mapsto x(t)\,$ of 
$\,v=\navp$, and set $\,(\,\,)\dot{\,}=\,d/dt\,$ (which is applied to 
functions restricted to the curve). Our initial task is to find conditions on 
$\,W=-\hs d\fc/d\vp\,$ and $\,Y\nnh=\Delta\vp$, restricted to the curve, 
necessary and sufficient for the linear ordinary differential equation
\begin{equation}\label{ode}
\dot\fy\,+\,\fy Y\,=\hs\,W
\end{equation}
to have a solution $\,\fy:\bbR\to\bbR\,$ with finite limits 
$\,\fy(\pm\infty)$. Such $\,\fy$, if it exists, must be unique. This is 
obvious from the preceding paragraph (which actually shows more: namely, there 
is at most one solution $\,\fy\,$ with a finite limit at $\,\infty$, and at 
most one with a finite limit at $\,-\infty$).

Writing $\,\int_a^{\hs b}[h]\,$ instead of $\,\int_a^{\hs b}h(t)\,dt\,$ 
whenever $\,h:\bbR\to\bbR\,$ is a continuous function and 
$\,a,b\in\bbR\cup\{\infty,-\infty\}$, we see that the condition
\begin{equation}\label{con}
\textstyle{\int_{-\infty}^{\hs\infty}}[\hh\zeta W\hh]=0\hs,\hskip5pt
\mathrm{where}\hskip4pt\zeta=e^Z\hskip4pt\mathrm{and}
\hskip4ptZ:\bbR\to\bbR\hskip5pt\mathrm{is\ an\ antiderivative\ 
of}\hskip4ptY,
\end{equation}
is necessary for the existence of a bounded solution $\,\fy\,$ to (\ref{ode}), 
as well as sufficient for (\ref{ode}) to have a solution $\,\fy\,$ with finite 
limits $\,\fy(\pm\infty)$. (Under our assumptions, $\,\zeta\,$ is always 
integrable, and hence so is $\,\zeta W$, while the equality in (\ref{con}) 
remains unchanged if one replaces $\,Z\,$ by another antiderivative of 
$\,Y\nnh$.) In fact, necessity of (\ref{con}) follows since, multiplying both 
sides of (\ref{ode}) by $\,\zeta$, we can rewrite (\ref{ode}) as 
$\,(\zeta\fy)\dot{\hskip2pt}\nh=\zeta W\nnh$, while boundedness of $\,\fy\,$ 
implies that $\,\zeta\fy\,$ has limits equal to $\,0\,$ at both $\,\infty\,$ 
and $\,-\infty$, so that integrating the last equality we get (\ref{con}). For 
sufficiency, note that, in view of l'Hospital's rule, 
$\,\fy(\pm\infty)=W(\pm\infty)/Y(\pm\infty)\,$ if one defines the solution 
$\,\fy\,$ to (\ref{ode}) by 
$\,\fy(t)=[\hs\zeta(t)]^{-1}\int_{-\infty}^{\hs t}[\hh\zeta W\hh]$. 

The requirement that (\ref{con}) hold along every maximal integral curve of 
$\,v$, as a restriction on the choice of a nonconstant $\,C^\infty$ function 
$\,\fc:[\hh\vp_{\mathrm{min}},\vp_{\mathrm{max}}]\to\bbR\,$ that would become 
the function $\,\fc\,$ of Lemma~\ref{spbic}(ii) for a nontrivial special 
bi\-con\-for\-mal change of $\,(g,\nnh\vp)$, is therefore necessary for such a 
bi\-con\-for\-mal change to exist. How restrictive this requirement is 
depends on $\,(g,\nnh\vp)$. For instance, if $\,\vp\,$ is a special 
K\"ah\-ler-Ric\-ci potential on $\,(M,g)$, (\ref{con}) states that 
$\,W\nh$, as a function of 
$\,\vp\in[\hh\vp_{\mathrm{min}},\vp_{\mathrm{max}}]$, should be 
$\,L^2$-or\-thog\-o\-nal to just one specific function of $\,\vp$. In 
general, however, the dependence of $\,\vp\,$ on $\,t\,$ varies with the 
integral curve, so that (\ref{con}) amounts to a much stronger 
$\,L^2$-or\-thog\-o\-nal\-i\-ty condition.

%\section{A nonexistence result}\label{ne}
%\setcounter{equation}{0}

\section{Remarks on the Ric\-ci form and scalar curvature}\label{rr}
\setcounter{equation}{0}
Let $\,g\,$ be a K\"ah\-ler metric on a complex manifold $\,M\,$ of complex 
dimension $\,m\ge2$. The Ric\-ci form of $\,g\,$ then is given by 
$\,\rho=\hh\mathrm{Ric}\hh(J\,\cdot\,,\,\cdot\,)$. The Ric\-ci forms of two 
K\"ah\-ler metrics $\,g,\hatg\,$ on $\,M\,$ are related by 
$\,\haro=\rho-i\hskip1pt\partial\overline{\partial}\,\log\gamma$, where 
$\,\hato\hs^{\wedge m}\nh=\gamma\hs\omega\hs^{\wedge m}\nnh$, that is,
$\,\gamma:M\to(0,\infty)\,$ is the ratio of the volume elements. If 
$\,(g,\nnh\vp)\,$ and $\,(\hatg,\havp)\,$ are met\-\hbox{ric-}\hskip.7ptpo\-ten\-tial pairs on 
$\,M\nh$, with a special bi\-con\-for\-mal change (\ref{bcm}), this yields
\begin{equation}\label{rcf}
\haro\,\,=\,\,\rho\hs\,
-\,\hs(m-1)\hs i\hskip1pt\partial\overline{\partial}\,\log\fx\hs\,
-\,\hs i\hskip1pt\partial\overline{\partial}\,\log(\fx-Q\hh\fy)\hs,
\end{equation}
since $\,\gamma=(\fx-Q\hh\fy)\hs\fx^{m-1}\nnh$, cf.\ (\ref{eig}).

When $\,m=2$, we have $\,\gamma=(\fx-Q\hh\fy)\hs\fx\nh$, as well as 
$\,4\hh\rho\wedge\hs\omega=\se\hs\omega\wedge\hs\omega$, while 
$\,4\hh(i\hskip1pt\partial\overline{\partial}\hs\psi)\wedge\hs\omega
=(\Delta\psi)\hs\omega\wedge\hs\omega\,$ for any function $\,\psi$, and 
$\,4\hh\rho\wedge\hs\xi\wedge\hh d\vp
=-(Q\hh\se\hh+\hh d_vY)\hs\omega\wedge\hs\omega$. The last three relations are 
direct consequences of the easily-verified formula
\begin{equation}\label{zwa}
4\hs\zeta\wedge\alpha\wedge\hs\xi\,\,
=\,\,[(\mathrm{tr}_\bbR^{\phantom i}A)\hh g(v,v)-2\hh g(\hn Av,v)]
\hs\omega\wedge\hs\omega
\end{equation}
valid whenever $\,(M,g)\,$ is a K\"ah\-ler surface, $\,\omega\,$ stands for 
its K\"ah\-ler form, $\,v\,$ is a tangent vector field, $\,A:\tm\to\tm\,$ is 
a bundle morphism commuting with the complex structure tensor $\,J\,$ and 
self-ad\-joint at every point, while $\,\zeta=g(JA\hs\cdot\,,\,\cdot\,)$, 
$\,\alpha=g(v,\,\cdot\,)\,$ and $\,\xi=g(Jv,\,\cdot\,)$. Note that 
(\ref{zwa}) gives $\,4\hs\zeta\wedge\hs\omega\hs
=(\mathrm{tr}_\bbR^{\phantom i}A)\hs\omega\wedge\hs\omega$. Hence, by 
(\ref{rcf}) with $\,m=2\,$ and (\ref{hoe}), the scalar curvatures $\,\se\,$ of 
$\,g\,$ and $\,\hasc\,$ of $\,\hatg\,$ satisfy the relation
\begin{equation}\label{hsc}
\gamma\hs\hasc\,\,=\,\,(\fx\,-\,\hh Q\hh\fy)\hs(\se\hs-\Delta\log\nh\gamma)\,
+\,\fy\hs d_v(d_v\log\nh\gamma-Y)\hs.
\end{equation}
The equalities $\,\gamma=(\fx-Q\hh\fy)\hs\fx\,$ and 
$\,\fx\nh-Q\hh\fy=\fc(\vp)\,$ (see Lemma~\ref{spbic}(ii)) make it possible to 
rewrite (\ref{hsc}) in a number of ways.

\setcounter{section}{1}
\renewcommand{\thesection}{\Alph{section}}
\section*{Appendix: Kil\-ling potentials with geodesic gradients}%\label{ae}
\setcounter{equation}{0}
\setcounter{theorem}{0}
The following construction generalizes that of 
\cite[\S5]{derdzinski-maschler-06} (in the case $\,m=2$), and gives rise to 
compact K\"ah\-ler surfaces $\,(M,g)\,$ with nonconstant Kil\-ling potentials 
$\,\vp$, which have geodesic gradients, but need {\it not\/} be special 
K\"ah\-ler-Ric\-ci potentials.

One begins by fixing a nontrivial closed interval 
$\,\bbI=[\hh\vp_{\mathrm{min}},\vp_{\mathrm{max}}]$, a constant 
$\,a\in(0,\infty)$, a compact K\"ah\-ler manifold $\,(N,h)\,$ of complex 
dimension $\,1\,$ (or, equivalently, a closed oriented real surface $\,N\hs$ 
endowed with a Riemannian metric $\,h$), $\,C^\infty$ mappings 
$\,\bbI\ni\vp\mapsto Q\in\bbR\,$ and 
$\,\yj:N\nh\to\bbRP^1\nh\smallsetminus\bbI\,$ such that $\,Q=0\,$ at the 
endpoints $\,\vp_{\mathrm{min}},\vp_{\mathrm{max}}$ and $\,Q>0\,$ on the open 
interval $\,\bbI^\circ\nnh=(\vp_{\mathrm{min}},\vp_{\mathrm{max}})$, while 
$\,Q\hh'\nh=2a\,$ at $\,\vp_{\mathrm{min}}$ and $\,Q\hh'\nh=-2a\,$ at 
$\,\vp_{\mathrm{max}}$. The use of the symbol $\,\yj\,$ conforms to the 
notations of \cite[\S5]{derdzinski-maschler-06}, where 
$\,\yj\in\bbR\smallsetminus\bbI
\subset\bbRP^1\nh\smallsetminus\bbI\,$ was a real constant. Here and 
below, $\,(\hskip2.5pt)'\nh=d/d\vp\,$ and $\,\bbR\,$ is treated as a subset of 
$\,\bbRP^1$ via the usual embedding $\,\vp\mapsto[\vp,1]\,$ (the brackets 
denoting, this time, the homogeneous coordinates in $\,\bbRP^1$). For 
algebraic operations involving $\,\infty=[\hh1,0\hs]\in\bbRP^1$ and elements 
of $\,\bbR\subset\bbRP^1\nnh$, the standard conventions apply; thus, 
$\,\infty^{-1}\nnh=0$. Since we need a canonically selected point $\,\vp_*$ in 
$\,\bbI$, we choose $\,\vp_*$ to be the midpoint of $\,\bbI$.

In addition, let us fix a $\,C^\infty$ complex line bundle $\,\mathcal{L}\,$ 
over $\,N\hs$ along with a Hermitian fibre metric $\,\langle\,,\rangle\,$ in 
$\,\mathcal{L}$, and a connection in $\,\mathcal{L}\,$ making 
$\,\langle\,,\rangle\,$ parallel and having the curvature form 
$\,\varOmega\hs=-\hs a\hs(\vp_*-\yj)^{-1}\omh\nh$, where $\,\omh$ is the 
K\"ahler form of $\,(N,h)$. (Thus, $\,\varOmega=0\,$ at points at which 
$\,\yj=\infty$.) The symbol $\,\mathcal{L}\,$ also denotes the total space of 
the bundle, while $\,\mathcal{V}\,$ and $\,\mathcal{H}\,$ stand for the 
vertical distribution $\,\mathrm{Ker}\hskip2.7ptd\proj\,$ and the horizontal 
distribution of our connection, $\,\proj\,$ being the projection 
$\,\mathcal{L}\to N\nh$. Treating the norm function 
$\,r:\mathcal{L}\to[\hs0,\infty)\,$ of $\,\langle\,,\rangle$, simultaneously, 
as an independent variable ranging over $\,[\hs0,\infty)$, we finally select a 
$\,C^\infty$ dif\-feo\-mor\-phism 
$\,\bbI^\circ\nh\ni\vp\mapsto r\in(0,\infty)\,$ such that 
$\,dr/d\vp=ar/Q$.

The above data allow us to define a Riemannian metric $\,g\,$ on 
$\,M'\nh=\mathcal{L}\smallsetminus N\nh$, where $\,N\hs$ is identified with 
the zero section, by 
$\,g=(\vp_*-\yj\circ\proj)^{-1}(\vp-\yj\circ\proj)\hs\proj^*\nh h\,$ or 
$\,g=\proj^*\nh h\,$ on $\,\mathcal{H}$, 
$\,g=(ar)^{-2}Q\,\text{\rm Re}\hskip1pt\langle\,,\rangle\,$ on 
$\,\mathcal{V}\nnh$, and $\,g(\mathcal{H},\mathcal{V})=\{0\}$. On 
$\,\mathcal{H}$, the first formula is to be used in the $\,\proj$-pre\-im\-age 
of the set in $\,N\hs$ on which $\,\yj\ne\infty$, and the second one on its 
complement. Note that, due to the fixed dif\-feo\-mor\-phic correspondence 
between the variables $\,\vp\,$ and $\,r\nh$, we may view $\,\vp\,$ (and hence 
$\,Q$) as a function $\,M'\nh\to\bbR$, while 
$\,C^\infty\nnh$-dif\-fer\-en\-tia\-bil\-i\-ty of the algebraic operations in 
$\,\bbRP^1\nnh$, wherever they are permitted, implies that $\,g\,$ is of class 
$\,C^\infty\nnh$.

The vertical vector field $\,v\,$ on $\,\mathcal{L}$, the restriction of which 
to each fibre of $\,\mathcal{L}\,$ equals $\,a\,$ times the radial (identity) 
vector field on the fibre, is easily seen to have the property that 
$\,d_v=Q\,d/d\vp$, with both sides viewed as operators acting on $\,C^\infty$ 
functions of $\,\vp$. Hence $\,v=\navp$, that is, $\,v\,$ is the 
$\,g$-gradient of $\,\vp$.

Clearly, $\,(M'\nnh,g)\,$ becomes an almost Her\-mit\-i\-an manifold when 
equipped with the unique almost complex structure $\,J\,$ such that the 
sub\-bun\-dles $\,\mathcal{V}\,$ and $\,\mathcal{H}\,$ of $\,\tm'$ are 
$\,J$-in\-var\-i\-ant and $\,J_x$ restricted to $\,\mathcal{V}\hskip-2.3pt_x$, 
or $\,\mathcal{H}_x$, for any $\,x\in M'\nnh$, coincides with the complex 
structure of the fibre $\,\mathcal{L}_{\proj(x)}$ or, respectively, with the 
$\,d\proj_x$-pull\-back of the complex structure of $\,N\nh$.

If $\,M\,$ now denotes the $\,\bbCP^1$ bundle over $\,N\hs$ obtained as the 
projective compactification of $\,\mathcal{L}$, then $\,g,\vp\,$ and $\,J\,$ 
have $\,C^\infty$ extensions to a metric, function and almost complex 
structure on $\,M\nh$, still denoted by $\,g,\vp\,$ and $\,J$. In addition, 
$\,g\,$ is a K\"ahler metric, that is, $\,\nabla J=0$, while $\,\vp\,$ is a 
Kil\-ling potential with a geodesic gradient on the compact K\"ah\-ler surface 
$\,(M,g)$, but, unless the function 
$\,\yj:N\nh\to\bbRP^1\nh\smallsetminus\bbI\,$ is constant, $\,\vp\,$ is not a 
special K\"ah\-ler-Ric\-ci potential. For details, see \cite{derdzinski-11}.

\ 

\noindent{\it Proof of Theorem~\ref{geodg}}\hskip4ptThe following 
classification theorem was established in \cite{derdzinski-11}:

{\it Let\/ $\,\vp\,$ be a nonconstant Kil\-ling potential with a geodesic 
gradient on a compact K\"ah\-ler surface\/ 
$\,(M,g)$. If\/ $\,\vp\,$ is not a special K\"ah\-ler-Ric\-ci potential on\/ 
$\,(M,g)$, then, up to a bi\-hol\-o\-mor\-phic isometry, the triple\/ 
$\,(M,g,\vp)\,$ arises from the above construction applied to some data\/ 
$\,\bbI,a,N\nh,h,\mathcal{L},\langle\,,\rangle,\mathcal{H},\yj\,$ and\/ 
$\,\vp\mapsto Q\,$ with the required properties, such that the function\/ 
$\,\yj:N\nh\to\bbRP^1\nh\smallsetminus\bbI\,$ is nonconstant.}

We may thus assume that $\,(M,g)\,$ and $\,\vp\,$ are the objects constructed 
above. For $\,\fc\,$ as in the statement of Theorem~\ref{geodg} and any fixed 
$\,y\in N\nh$, let 
$\,F\nnh_y:[\hh\vp_{\mathrm{min}},\vp_{\mathrm{max}}]\to\bbR\,$ be the 
antiderivative, vanishing at $\,\vp_{\mathrm{min}}$, of the function 
$\,-\hs[\vp-\yj(y)]\hh\fc'(\vp)\,$ of the variable $\,\vp$. Thus, 
$\,F\nnh_y=0\,$ at both endpoints $\,\vp_{\mathrm{min}},\vp_{\mathrm{max}}$. 
Due to the boundary conditions imposed on $\,Q\,$ and the first-or\-der Taylor 
formula, $\,F\nnh_y$ is is smoothly divisible by $\,Q\,$ on the whole closed 
interval $\,[\hh\vp_{\mathrm{min}},\vp_{\mathrm{max}}]$, that is, 
$\,QE_y(\vp)=F\nnh_y(\vp)\,$ for some $\,C^\infty$ function $\,E_y$, and we 
may define a $\,C^\infty$ function $\,\fy:M\to\bbR\,$ by 
$\,\fy(x)=E_y(\vp)/[\vp-\yj(y)]$. (Here $\,\vp\,$ stands for $\,\vp(x)$, and 
$\,y=\proj(x)$, with  $\,\proj:M\to N\hs$ denoting the bundle projection.)

Next, $\,Y\nnh=\Delta\vp\,$ is given by 
$\,Y\nnh=(\vp-\yj\circ\proj)^{-1}Q+\hs dQ/d\vp\,$ (see \cite{derdzinski-11}). 
Since, as we noted above, $\,d_v=Q\,d/d\vp$, condition (iii) of 
Lemma~\ref{spbic} follows. Adding a constant to $\,\fc\nh$, we also obtain 
(iv) in Lemma~\ref{spbic}, if $\,P\nh,\fx\,$ and $\,\havp\,$ chosen so as to 
satisfy (ii) and (i) in Lemma~\ref{spbic}. This completes the proof.
\qed

\begin{acknowledgements}
The author thanks Gideon Maschler for helpful comments and suggestions, 
including the italicized observation at the end of Section~\ref{sk}.
\end{acknowledgements}


\begin{thebibliography}{19} % '2nd argument contains the widest acronym'

\bibitem{calabi}Calabi, E.: Extremal K\"ah\-ler metrics. In: Yau, 
S.-T. (ed.) Seminar on Differential Geometry, pp.\ 259--290. Annals of Math.\ 
Studies \textbf{102}, Princeton Univ.\ Press, Princeton, NJ (1982) 

\bibitem{cao}Cao, H.-D.: Existence of gradient K\"ah\-ler-Ric\-ci 
sol\-i\-tons. In: Chow, B., et al. (eds.) Elliptic and Parabolic Methods in 
Geometry, Minneapolis, MN, 1994, pp.\ 1--16. A.K. Peters, Wellesley, MA (1996)

\bibitem{chen-lebrun-weber}Chen, X., LeBrun, C., 
Weber, B.: On con\-for\-mal\-ly K\"ah\-ler, Einstein manifolds. 
J.\ Amer.\ Math.\ Soc. \textbf{21}, 1137--1168 (2008)

\bibitem{chow}Chow, B.: Ric\-ci flow and Ein\-stein metrics 
in low dimensions. In: LeBrun, C., Wang, M. (eds.) Surveys in Differential 
Geometry, Vol.VI, pp.\  187--220. Internat.\ Press, Boston, MA 
(1999)

\bibitem{derdzinski-83}Derdzi\'nski, A.: Self-du\-al K\"ah\-ler 
manifolds and Ein\-stein manifolds of dimension four. Compos.\ Math. 
\textbf{49}, 405--433 (1983)

\bibitem{derdzinski-11}Derdzinski, A.: Kil\-ling potentials with geodesic 
gradients, in preparation, current version available from 
http:/\hskip-.5pt/www.math.ohio-state.edu/\hs\~{}andrzej\hs/gg.pdf

\bibitem{derdzinski-maschler-03}Derdzinski, A., Maschler, G.: Local 
classification of con\-for\-mal\-ly-Ein\-stein K\"ah\-ler metrics in higher 
dimensions. Proc.\ London Math.\ Soc. \textbf{87}(3), 779--819 (2003)

\bibitem{derdzinski-maschler-06}Derdzinski, A., Maschler, G.: Special 
K\"ah\-ler-Ric\-ci potentials on compact K\"ah\-ler manifolds. J.\ reine 
angew.\ Math. \textbf{593}, 73--116 (2006)

\bibitem{el-mansouri}El-Mansouri, J.: On adapted 
bi-con\-for\-mal metrics. J.\ Geom. \textbf{81}, 30--45 (2004)

\bibitem{ganchev-mihova}Gan\-chev, G., Mi\-ho\-va, V\nnh.: 
K\"ah\-ler manifolds of qu\-asi-con\-stant hol\-o\-mor\-phic sectional 
curvatures. Cent.\ Eur.\ J.\ Math. \textbf{6}, 43--75 (2008)

\bibitem{hamilton}Hamilton, R.S.: The Ric\-ci flow 
on surfaces. In: Isenberg, J.A. (ed.) Mathematics and General Relativity, 
Santa Cruz, CA, 1986, pp.\ 237--262. Contemp.\ Math., vol.\ \textbf{71}, 
Amer.\ Math.\ Soc., Providence, RI (1988)

\bibitem{jelonek}Jelonek, W.: K\"a\-hler manifolds with 
qua\-si-con\-stant hol\-o\-mor\-phic curvature. Ann.\ Global Anal.\ Geom. 
\textbf{36}, 143--159 (2009)

\bibitem{kobayashi-nomizu}Kobayashi, S., Nomizu, K.: 
Foundations of Differential Geometry, Vol.\ I. Interscience, New York (1963)

\bibitem{koiso}Koiso, N.: On rotationally symmetric Hamilton's
equation for K\"ah\-ler-Ein\-stein metrics. In: Ochiai, T. (ed.) Recent Topics 
in Differential and Analytic Geometry, pp.\ , 327--337. Adv.\ Stud.\ Pure 
Math., vol.\ \textbf{18}-{\rm I}, Academic Press, Boston, MA (1990)

\bibitem{lebrun}LeBrun, C.: Einstein metrics on complex surfaces. 
In: Andersen, J.E., et al. (eds.) Geometry and Physics, Aarhus, 1995, pp.\ 
167--176. Lecture Notes in Pure and Applied Mathematics \textbf{184}, Dekker, 
New York (1997) 

\bibitem{matsushima}Matsushima, Y.: Sur la structure du groupe
d'hom\'eomorphismes analytiques d'une certaine vari\'et\'e kaehl\'erienne. 
Nagoya Math.\ J. \textbf{11}, 145-150 (1957)

\bibitem{nicolaescu}Nicolaescu, Liviu I.: An Invitation to Morse Theory. 
Universitext, Springer, New York (2007)

\bibitem{page}Page, D.: A compact rotating gravitational in\-stant\-on. 
Phys.\ Lett.\textbf{79B}, 235--238 (1978)

\bibitem{perelman}Perelman, G.: The entropy formula for the 
Ric\-ci flow and its geometric applications. Preprint, %available from 
arXiv:math.DG/0211159.

\bibitem{tian-zhu}Tian, G., Zhu, X.: A new hol\-o\-mor\-phic invariant and 
uniqueness of K\"ah\-ler-Ric\-ci solitons. Comment.\ Math.\ Helv. \textbf{77}, 
297--325 (2002)

\bibitem{wang-zhu}Wang, X.J., Zhu, X.: K\"ah\-ler-Ric\-ci solitons 
on toric manifolds with positive first Chern class. Adv.\ in Math. 
\textbf{188}, 87--103 (2004)

\end{thebibliography}
\end{document}